\documentclass{elsarticle}

\usepackage{morefloats}
\usepackage{graphicx}
\usepackage{epsfig}
\usepackage{import}
\usepackage{psfrag}
\usepackage{enumerate}
\usepackage{color}
\usepackage{amsmath}
\usepackage{amsfonts}
\usepackage{amssymb}
\usepackage{mathrsfs}
\newcommand\calE{\mathscr E}
\newcommand\calR{\mathscr R}

\newcommand\calF{\mathscr F}
\newcommand\calM{\mathscr M}
\newcommand\N{\mathbb N}
\newcommand\R{\mathbb R}
\newcommand\pl{\partial}
\renewcommand\d{\mathrm d}
\usepackage{minibox}
\newcommand\JUMP[2]{\mathchoice
                  {\big[\hspace*{-.3em}\big[#1\big]\hspace*{-.3em}\big]_{#2}}
                   {[\hspace*{-.15em}[#1]\hspace*{-.15em}]_{#2}}
                   {[\![#1]\!]_{#2}}
                   {[\![#1]\!]_{#2}}}
\newcommand\JUMPPOWER[3]{\mathchoice
                  {\big[\hspace*{-.3em}\big[#1\big]\hspace*{-.3em}\big]_{#2}^{#3}}
                   {[\hspace*{-.15em}[#1]\hspace*{-.15em}]_{#2}^{#3}}
                   {[\![#1]\!]_{#2}^{#3}}
                   {[\![#1]\!]_{#2}^{#3}}}
\newcommand{\GD}{\mathchoice
                  {\Gamma_{\hspace*{-.15em}\mbox{\tiny\rm D}}}
                  {\Gamma_{\hspace*{-.15em}\mbox{\tiny\rm D}}}
                  {\Gamma_{\hspace*{-.1em}\mbox{\tiny\rm D}}}
                  {\Gamma_{\hspace*{-.05em}\mbox{\tiny\rm D}}}}
\newcommand{\GN}{\mathchoice
                  {\Gamma_{\hspace*{-.15em}\mbox{\tiny\rm N}}}
                  {\Gamma_{\hspace*{-.15em}\mbox{\tiny\rm N}}}
                  {\Gamma_{\hspace*{-.1em}\mbox{\tiny\rm N}}}
                  {\Gamma_{\hspace*{-.05em}\mbox{\tiny\rm N}}}}
\newcommand\GC{\mathchoice{\Gamma_{\hspace*{-.15em}\mbox{\tiny\rm C}}}
                          {\Gamma_{\hspace*{-.15em}\mbox{\tiny\rm C}}}
                          {\Gamma_{\hspace*{-.05em}\mbox{\tiny\rm C}}}
                          {\Gamma_{\hspace*{-.05em}\mbox{\tiny\rm C}}}}
\newcommand{\GDbar}{\mathchoice
                  {\overline{\Gamma}_{\hspace*{-.15em}\mbox{\tiny\rm D}}}
                  {\overline{\Gamma}_{\hspace*{-.15em}\mbox{\tiny\rm D}}}
                  {\overline{\Gamma}_{\hspace*{-.1em}\mbox{\tiny\rm D}}}
                  {\overline{\Gamma}_{\hspace*{-.05em}\mbox{\tiny\rm D}}}}
\newcommand\GCbar{\mathchoice{\overline{\Gamma}_{\hspace*{-.15em}\mbox{\tiny\rm C}}}
                          {\overline{\Gamma}_{\hspace*{-.15em}\mbox{\tiny\rm C}}}
                          {\overline{\Gamma}_{\hspace*{-.05em}\mbox{\tiny\rm C}}}
                          {\overline{\Gamma}_{\hspace*{-.05em}\mbox{\tiny\rm C}}}}
\newcommand\DT[1]{\mathchoice
                 {{\buildrel{\hspace*{.1em}\text{\LARGE.}}\over{#1}}}
                 {{\buildrel{\hspace*{.1em}\text{\Large.}}\over{#1}}}
                 {{\buildrel{\hspace*{.1em}\text{\large.}}\over{#1}}}
                 {{\buildrel{\hspace*{.1em}\text{\large.}}\over{#1}}}}
\newcommand{\Item}[2]{\parbox[t]{.05\textwidth}{#1}\hfill%
      \parbox[t]{.95\textwidth}{#2}\vspace*{.8mm}}
\newcommand{\DDD}[3]{\begin{array}[t]{c}#1\vspace*{-1em}\\_{#2}\vspace*{-.3em}\\_{#3}\end{array}}
\newcommand{\ddd}[3]{\DDD{\begin{array}[t]{c}\underbrace{#1}\vspace*{.6em}\end{array}}{\text{\footnotesize #2}}{\text{\footnotesize #3}}}
\newcommand\Diss{\mathrm{Diss}}
\newcommand\uD{u_{\text{\tiny\rm D}}}
\newcommand\pN{p_{\text{\tiny\rm N}}}
\newcommand\uC{u_{\text{\tiny\rm C}}}

\newcommand\uDi{u^i_{\text{\tiny\rm D}}}
\newcommand\pNi{p^i_{\text{\tiny\rm N}}}

\newcommand{\GDi}{\Gamma^i_{\tiny\rm D}}
\newcommand{\GNi}{\Gamma^i_{\tiny\rm N}}

\newcommand\duD{\DT{u}_{\text{\tiny\rm D}}}
\newcommand\dpN{\DT{p}_{\text{\tiny\rm N}}}

\newcommand\pD{p_{\text{\tiny\rm D}}}
\newcommand\uN{u_{\text{\tiny\rm N}}}
\newcommand\pC{p_{\text{\tiny\rm C}}}

\newcommand\SD{P^{\text{\tiny\rm D}}}
\newcommand\SO{P}
\newcommand\SN{P^{\text{\tiny\rm N}}}
\newcommand\SC{P^{\text{\tiny\rm C}}}

\newcommand\GIC{G_{_{\rm I \tiny c}}}
\newcommand\GIIC{G_{_{\rm II \tiny c}}}

\newcommand\COL[1]{{#1}}

\newcommand\uCS{u^{\text{\tiny\rm C}}}
\newcommand\uDS{u^{\text{\tiny\rm D}}}
\newcommand\uNS{u^{\text{\tiny\rm N}}}

\journal{Computational Mechanics}

\begin{document}

\begin{frontmatter}


 \author{C.G.~Panagiotopoulos\corref{cor2}\fnref{label1}}
 \ead{cpanagiotopoulos@us.es}
 
 \author{V.~Manti\v c\corref{cor1}\fnref{label1}}
 \ead{mantic@etsi.us.es}
 \address{Group of Elasticity and Strength of Materials, Department of Continuum Mechanics, School of Engineering, University of Seville, Camino de los Descubrimientos s/n \\ES-410 92 Sevilla, Spain}

\author{T.~Roub\'\i\v cek\corref{cor3}\fnref{label2}}
\ead{roubicek@karlin.mff.cuni.cz}
 \address{Mathematical Institute, Charles University, Sokolovsk\'a 83, CZ-186~75~Praha~8, Czech Republic}
 \address{Institute of Thermomechanics of the ASCR, Dolej\v skova 5, CZ-182 00 Praha 8, Czech Republic}

\fntext[label1]{The support by the Junta de Andaluc\'{\i}a
and Fondo Social Europeo (Proyecto de Excelencia TEP-4051)
is warmly acknowledged. VM also acknowledges the support by the Ministerio de Ciencia e Innovaci\'on
(Proyecto MAT2009-14022).} 
\fntext[label2]{TR~acknowledges partial support from the grants 201/09/0917, 201/10/0357, and 201/12/0671 (GA \v CR), and the institutional support RVO: 67985971 (\v CR).}
 
\title{BEM solution of delamination problems using an interface damage and plasticity model}
\fntext[label3]{DOI:10.1007/s00466-012-0826-3}

\begin{abstract}\small{
The problem of quasistatic and rate-in\-de\-pend\-ent evolution of elastic-plastic-brittle delamination at small strains is considered. Delamination processes for linear elastic bodies glued by an adhesive to each other or to a rigid outer surface are studied. The energy amounts dissipated in fracture Mode I (opening) and Mode II (shear) at an interface  may be different. A concept of internal parameters is used here on the delaminating interfaces, involving a couple of scalar damage variable and a plastic tangential slip with kinematic-type hardening. The so-called energetic solution concept is employed.
An inelastic process at an interface is devised in such a way that the dissipated energy depends only on the rates of internal parameters and therefore the model is associative. A fully implicit time discretization is combined with a spatial  discretization of elastic bodies by the BEM to solve the delamination problem.
The BEM is used in the solution of the respective boundary value problems, for each subdomain separately, to compute the corresponding total potential energy. Sample problems are analysed by a collocation BEM code to illustrate the capabilities of the numerical procedure developed.}
\end{abstract}

\begin{keyword}
\small{Interface fracture} \sep 
\small{rate-independent quasistatic model} \sep 
\small{adhesive contact} \sep 
\small{energetic solutions} \sep 
\small{delamination/debonding} \sep 
\small{plastic slip} \sep 
\small{boundary element method}
\end{keyword}

\end{frontmatter}


\section{Introduction} \label{intro}
Applications of layered structures are numerous and continuously increased, an example being the massive use of composite materials in aeronautical industry at present. Usually the interfaces
between these rather bulk laminates consist of very thin adhesive layers. For efficient computations,
these adhesive layers may be approximated by zero thickness interface layers.
There are many situations where an adhesive layer is found to be partially or fully damaged. This process is frequently referred to as \emph{delamination} or \emph{debonding} of adjacent material laminas.
In this work the description of the damage is based on a scalar \emph{damage} quantity (variable, cf.\ \cite{Frem85DAS}), which is defined at   interfaces   and takes values from the interval $[0,1]$, with zero value meaning no adhesion due to the total damage of the adhesive while the
unit value  meaning complete operation of the adhesive without any damage. During a damage   evolution   the damage variable decays in time, and it is assumed that a specific amount of energy has to be released (dissipated). This simplified approach, motivated essentially by Griffith \cite{Grif20PRFS}, is often inadequate as it is observed experimentally that considerably more energy is usually needed to perform delamination in shear Mode II than in opening Mode I \cite{BanAsk00NFCI,HutSuo92MMCL,LieCha92ASIF,TveHut93IPMM}. Motivated by the microscopical idea of interface plasticity \cite{LieCha92ASIF,TveHut93IPMM}, an extra inelastic parameter is  introduced \cite{RoKrZe11DACM,RoMaPa??QMMD}, which describes some plastic slip that may occur in the tangent direction of an interface before its debonding.

An alternative approach to model fracture-mode-sensitive
delamination uses only the delamination variable but makes the
dissipated energy directly dependent on the so-called fracture mode mixity angle,
cf.~\eqref{FMMAE}--\eqref{FMMADS} below.
This approach has frequently been used in engineering models
\cite{TMGCP10ACTA,TMGP11BEMA}
but, it does not seem amenable to a rigorous mathematical analysis. In the present work we consider the delamination as a unidirectional process, i.e.\
no healing (or reconstruction) of adhesive is allowed, which covers
most of engineering applications.

The goal of this article is to  present and analyse from an
engineering  as well as numerical implementation  viewpoint some basic
features of the delamination model  devised in
\cite{RoKrZe11DACM,RoMaPa??QMMD} with different dissipated energies in Modes I
and   II.
In particular, in Section~\ref{sec:Theory} we briefly present the energetic approach employed. In Section~\ref{sec:model}, we concisely introduce the present model, while some engineering insight on this model is provided. Then, in Section~\ref{sec:numimpl}, the numerical implementation of the model, is presented.
Finally, in Section~\ref{sec:numexmp}, two-dimensional
simulations are developed, showing that the model is suitable for solving realistic  problems of delamination between  elastic  layers.

\section{Theoretical background}\label{sec:Theory}
\subsection{Problem definition}\label{sec:definition}
\begin{figure}
   \centering
   \def\svgwidth{.5\columnwidth}
   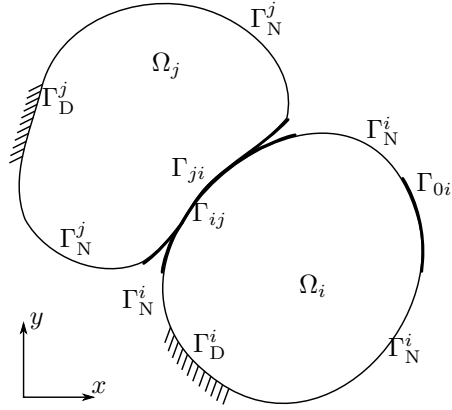
   \caption{Schematic illustration of the geometry and notation for a
two-dimensional case of two  bonded subdomains,
i.e.~$N{=}2$.}
\label{fig:geom}
\end{figure}

Let us consider an assemblage of $N$ elastic bodies, each of them defined
by a reference domain $\Omega_i$  $(i=1,...,N)$, with the Lipschitz boundary $\Gamma_i=\partial\Omega_i$, see Fig.~\ref{fig:geom}.
We denote by $\Gamma_{ij}=\partial\Omega_i \cap\partial\Omega_j$ the
(possibly empty) interface
boundary between $\Omega_i$ and $\Omega_j$  $(i,j=1,\ldots,N)$,
which may undergo delamination.
We also consider  possible delamination
on some parts of the outer boundary $\Gamma_{0i}$, which is assumed to be in Signorini elastic  contact with a fixed rigid surface, see Fig.~\ref{fig:geom}. The union of these parts is denoted as $\Gamma_{0}=\bigcup_{1\le i \le N} \Gamma_{0i}$.
We will denote $\GC:=\bigcup_{1\le i<j\le N}\Gamma_{ij}\cup\Gamma_{0}$.
 We assume that the rest of the outer boundary $\partial\Omega$ is the union of two disjoint subsets $\GD$ and $\GN$, where Dirichlet
(prescribed displacements $u_{\rm D}=u_{\rm D}(t)$)
and Neumann boundary conditions (prescribed tractions $p_{\rm N}=p_{\rm N}(t)$) are imposed, respectively. For the sake of simplicity of the following considerations, vanishing tractions $p_{\rm N}=0$ will be considered hereinafter, except for Sections \ref{subsec:BEM}, \ref{subsec:BoundFormPotEnergy} and \ref{par:TRACTIOMLOAD}.  The intersection of the closures of $\GC$ and $\GD$  is assumed to be the empty set, i.e. $\GCbar \cap \GDbar{=}\emptyset$. Any $\Gamma_{ij}$ is considered as an infinitely thin adhesive layer, represented by    springs  distributed continuously, similarly to the Winkler spring model, with distinct normal and tangential elastic stiffnesses of  values ranging from zero to infinity. Both the elastic subdomains and the adhesive layers are assumed to store energy, which is  given by a \emph{stored energy} functional $\calE(t,u,z)$ a function of time $t$, the displacements $u$ and the \emph{inelastic  (internal) parameters} collected in $z$. It is considered that two elastic subdomains $\Omega_i$ and $\Omega_j$, may debond along the interface $\Gamma_{ij}$. During this process the material of the adhesive can be damaged and plastified. The onset and growth of the damage and plastification, represented by the $z$ variables, does not depend on some internal time scale and therefore the process is considered as \emph{rate-independent}. The damage and plastification of the adhesive layer are accompanied by a release of stored energy. The \emph{dissipation potential} $\calR(\DT{z})$, with $\DT{z}:=\frac{\mathrm d z}{\mathrm d t}$, for a rate-independent process can be represented by a degree-1 homogeneous functional \cite{MiRoZe10CDEV}. The processes  described in this work  are assumed to be \emph{quasistatic}, i.e.   no inertia effects are taken into account.
The rate-independent evolution we have in mind is governed by the following
system of \emph{doubly nonlinear} degenerate abstract
\emph{static/evolution inclusions}, referred sometimes as Biot's equations
generalizing the original work \cite{Biot56TIT,Biot65MID}:
\begin{align}\label{Gm}
\pl_u\calE(t,u,z)\ni0\quad\text{and}\quad
\pl\calR\big(\DT{z}\big)+\pl_z\calE(t,u,z)
\ni0,
\end{align}
where the symbol
``$\partial$'' refers to a (partial) subdifferential,
relying on that $\calR(\cdot)$, $\calE(t,\cdot,z)$, and
$\calE(t,u,\cdot)$ are convex functionals. The first optimality condition of Eq.\,\eqref{Gm} represents the \emph{minimum energy principle}, while the latter one, the \emph{minimum dissipation potential principle}~\cite{RoKrZe11DACM}.

For the sake of simplicity, throughout this work, we will restrict ourselves
to the two-dimensional case, i.e.\  $\Omega_i\subset\R^2$ will be planar
domains, $i=1,...,N$, and
 $\Gamma_{ij}$ will be one-dimensional surfaces.

\subsection{Energetic solutions}
A fruitful concept of a certain weak solution to the doubly nonlinear
inclusion with degree-1 homogeneous dissipation potential $\calR$,
called energetic solutions, was developed by Mielke et al.\
\cite{MieThe04RIHM,MiThLe02VFRI}.
In the convex case, this concept is essentially equivalent to
conventional weak-solution concept, while in our case where
$\calE(t,\cdot,\cdot)$ is non-convex this concept represents
a certain generalization; cf.\ \cite{Miel05ERIS} for a
survey on the concept of energetic solutions and \cite{Miel09?DEMF}
for comparison with other concepts.

The process $(u({t}),z({t})),\ {t}\in[0,T]$ is
called an energetic solution to the initial-value problem \eqref{Gm}, if it
satisfies the following three conditions:
\\[0.5em]
\Item{(i)}{The \emph{energy equality}:}
\vspace*{-2em}
\begin{align}
&\ \ \ddd{\calE(T,u(T),z(T))_{_{_{_{_{_{}}}}}}\!\!\!}{stored energy}{at time $t=T$}
+ \!\!\ddd{\Diss_{\mathcal R}(z;[0,T])_{_{_{_{_{_{}}}}}}\!\!\!}{energy dissipated}{during $[0,T]$}
\label{total-energy}
\nonumber \\
 &\qquad\qquad\qquad =\!\!\ddd{\int_0^T\!\!\!\calE_t'(t,u,z)\,\d t}{work done by}{mechanical load}\!\!
+\!\!\ddd{\calE(0,u_0,z_0)_{_{_{_{_{_{_{}}}}}}}\!\!\!}{stored energy}{at time $t=0$}\!,
\end{align}
\Item{}{where}
\vspace*{-3em}
\begin{align}\label{def-of-diss}
\ \ \Diss_{\mathcal R}(z;[0,T]):=\sup
\sum_{j=1}^N\mathcal R(z(t_{j})-z(t_{j-1})),
\end{align}
\Item{}{with the supremum taken over all partitions\\
$0\le t_0<t_1<...<t_{N-1}\le t_N\le T$.
}
\Item{(ii)}{\emph{Stability inequality} for
any $t\in[0,T]$:}
\begin{align}\label{stability}
\ \ \calE(t,u,z)\le\calE(t,\tilde u,\tilde z)+\calR(\tilde z{-}z)
\ \ \ \text{ for any
}(\tilde u,\tilde z),
\end{align}
\Item{(iii)}{The \emph{initial conditions}:
$u(0)=u_0$ and $z(0)=z_0$.}
In Eq.\,\eqref{total-energy}, $\calE_t'$ is the partial derivative of $\calE$ with respect to time $t$.

\section{Model of interface damage and plasticity} \label{sec:model}
In this section we present the specific model adopted, in order to simulate
the nonlinear
inelastic behaviour of an adhesive layer, by defining a suitable
stored energy functional as well as a dissipation potential.
The present plastic-type model with  kinematic-type hardening \cite{HanRed99PMTN,simo-hughes} for
the
delamination problem, was devised essentially
in \cite{RoKrZe11DACM} without any mathematical or computational justification,
and further scrutinized in \cite{RoMaPa??QMMD}. Beside the
displacement $u$, two internal parameters are
used in order to describe the nonlinear behaviour of the adhesive:
the damage variable $\zeta$ and the plastic tangential slip variable $\pi$, which
together constitute the pair of inelastic variables $z=(\zeta,\pi)$.

\subsection{Stored energy}
Stored energy $\calE$ includes the elastic bulk contribution and the additional adhesive-surface contribution:
\begin{equation}\label{StoredEnergyDef}
{\calE}(t,u,z)={\calE}_{\rm el}(t,u)+{\calE}_{\rm adh}(u,z)
\end{equation}
 with
\begin{equation}\label{Eell}
{\calE}_{\rm el}(t,u)=\!\left\{\begin{array}{ll}
\sum_{i=1}^N\int_{\Omega_i}\mathbb C_i e(u){:}e(u)\,\d x,
\\
&\!\!\!\!\!\!\!\!\!\!\!\!\!\!\!\!\text{if }u|_{\GD}=u_{\rm D}(t),
\\[.3em]
   \infty&
\!\!\!\!\!\!\!\!\!\!\!\!\!\!\!\!\text{elsewhere.}
\end{array}\right.
\end{equation}
where $\mathbb C_i$ is the elastic moduli tensor in $\Omega_i$, and
\begin{align}
&\label{EnergyAdh}
{\calE}_{\rm adh}(u,z):=\!\left\{\begin{array}{ll}
  \displaystyle{\
\!\int_{\GC}\!\!
 \Big(\zeta\Big(
\frac{\kappa_{\rm n}}{2} \JUMPPOWER{u}{\rm n}{2}\!+
\frac{\kappa_{\rm t}}{2}\big(\JUMP{u}{\rm t}\!{-}\pi\big)^2\Big)}&\\[.9em]
\displaystyle{\ \ \ +\,\frac{\kappa_{_{\rm H}}}{2} \pi^2
+\frac{\kappa_0}r\big|\partial_{s}\zeta\big|^r  \Big)}
\displaystyle{\,\d S}
\\
&\!\!\!\!\!\!\!\!\!\!\!\!\!\!\!\!\!\!\!\!\!\!\!\!\!\!\!\!\!\!\!\!\!\!\!\!\!\!\!\!\!\!\!\!\!\!\!\!\!\text{if }0\le\zeta\le1\text{ and }\\
&\!\!\!\!\!\!\!\!\!\!\!\!\!\!\!\!\!\!\!\!\!\!\!\!\!\!\!\!\!\!\!\!\!\!\!\!\!\!\!\!\!\!\!\!\!\!\!\!\!\JUMP{u}{\rm n}\ge0\text{ on }\GC,
\\[.3em]
   \infty&\!\!\!\!\!\!\!\!\!\!\!\!\!\!\!\!\!\!\!\!\!\!\!\!\!\!\!\!\!\!\!\!\!\!\!\!\!\!\!\!\!\!\!\!\!\!\!\!\!\text{elsewhere,}
\end{array}\right.\hspace*{-1em}
\end{align}
where $\kappa_{\rm n}>0$ and  $\kappa_{\rm t}>0$ are the
 phenomenological elastic constants describing the stiffnesses of the linearly elastically
responding adhesive in the normal and tangential directions, respectively,
$\kappa_0>0$ is the so-called factor of influence of damage
\cite{MieRou06RIDP}, $\JUMP{u}{}=\JUMP{u}{\rm n}\nu+\JUMP{u}{\rm t}\tau$ with
$\JUMP{u}{\rm n}=\JUMP{u}{}{\cdot}\nu$ and $\JUMP{u}{\rm t}=\JUMP{u}{}{\cdot}\tau$, $\nu$ and $\tau$ being  unit normal and tangential vectors to $\GC$, and $\partial_{s}$ is  the tangential derivative   defined
on $\GC$).
For $\Gamma_{0}$ the outward normal $\nu$ is typically taken.
Constant parameter $\kappa_{_{\rm H}}$ stands for the  plastic modulus of kinematic hardening. Here we used the notation $\JUMP{u}{}$ for the differences of displacements
from both sides of $\GC$. We also assume $r>1$.
The last term
 in Eq.\,\eqref{EnergyAdh}, although bearing a physical interpretation
\cite{BazJir02NIFP}, is here introduced mainly for mathematical reasons
in order to facilitate
a proof of convergence; for further details see
also \cite{RoMaPa??QMMD},
but in specific simulations one may expect
reasonable numerical results
even if this term is neglected by setting $\kappa_0=0$.
The constraint $\JUMP{u}{\rm n}\ge0$ in Eq.\,\eqref{EnergyAdh} is actually the Signorini non-penetration condition of unilateral contact
\cite{KoMiRo06RIAD}.



\subsection{Dissipation potential}\label{subsect-dissip}
A suitable functional for the dissipation potential, describing both inelastic
processes of damage and plastic slip in the adhesive layer, and actually being a degree-1
homogeneous functional, is defined as:
\begin{align}\label{DissipationAdh}
&\calR(\DT{z})=\calR(\DT\zeta,\DT\pi):=
\begin{cases}\displaystyle{\int_{\GC}\!\!\GIC\,\big|\DT\zeta\big|
+\sigma_{\mathrm{t,yield}}\big|\DT\pi\big|\d S}
\\
&\!\!\!\!\!\!\!\!\!\!\!\!\!\!\!\!\!\!\!\!\!\!\!\!\!\!\!\!\!\!\!\!\!\!\!\!\!\!\!\!\!\!\!\!\!\!\text{ if }\DT\zeta\le0\text{ a.e.~on }\GC,\\
\infty&\!\!\!\!\!\!\!\!\!\!\!\!\!\!\!\!\!\!\!\!\!\!\!\!\!\!\!\!\!\!\!\!\!\!\!\!\!\!\!\!\!\!\!\!\!\!\text{ otherwise}.
\end{cases}
\end{align}
Parameter $\GIC>0$ is the
minimal energy required
for complete damage (debonding) of a unit area of the interface. In particular, we assume
it represents the interface fracture energy in Mode I.
Parameter  $\sigma_{\mathrm{t,yield}}>0$ is the interface yield
shear stress for initiation of tangential plastic slip along the interface.
The constraint $\DT\zeta\le0$ in \eqref{DissipationAdh}
makes the evolution of $\zeta$ irreversible, i.e.\ the model does not permit healing, which means that a debond appeared at some point can not be restored.

Note that, except trivial case when $u_{\rm D}$ is constant in time,
$\calE_t'$ in \eqref{total-energy} would not be well
defined. One way how to avoid this drawback, well consistent with BEM,
is to restrict the displacement only on $\GC$, assuming that
$\GC$ and $\GD$ are not touching each other. The restricted displacement
$u|_{\GC}$ will be denoted by $\uC$; in fact, in the case of
$\Gamma_{ij}$, it is a couple of traces
of $u$ from both sides of $\Gamma_{ij}$. As $u$ does not
occur in Eq.\,\eqref{DissipationAdh} and thus it is fully nondissipative,
Eq.\,\eqref{stability} implies that $u$ minimizes $\calE(t,\cdot,z)$
and thus, in fact, $\uC$ and $z$ determines $u$ at a given time $t$.
Thus, $\calE$ can be considered as a function of $\uC$ instead of $u$,
which makes $\calE_t'(t,\uC,z)$ well defined if $\uD$ is smooth in
time. On the other hand, we will not distinguish between $\JUMP{u}{}$ and
$\JUMP{\uC}{}$.
We will use this convention through the rest of this article.
\subsection{Engineering analysis of the traction-relative displacement law}
In the case of a linear elastic-brittle interface model \cite{TMGCP10ACTA,TMGP11BEMA}, the interface failure criterion  is connected to the \emph{energy release rate}
(ERR) concept.
It can be shown \cite{CORCARPU2009,LENCI2001} that the energy stored in the adhesive at the crack tip equals the ERR of a mixed mode crack propagating along a linear elastic interface, and can be evaluated
as:
\begin{equation}\label{ERR}
G=G_{\rm I}+G_{\rm II}=\frac{\kappa_{\rm n}\JUMPPOWER{u}{\rm n}{2}}{2}+\frac{\kappa_{\rm t}\JUMPPOWER{u}{\rm t}{2}}{2}.
\end{equation}
The so-called \emph{fracture mode mixity angles}, denoted as $\psi_G$, $\psi_u$, or $\psi_\sigma$, can be defined in terms of ERR as,
\begin{align}\label{FMMAE}
\tan^2{\psi_G}=\frac{G_{\rm II}}{G_{\rm I}},
\end{align}
as well as in terms of relative displacements and tractions, respectively,
\begin{align}\label{FMMADS}
\tan{\psi_u}=\frac{\JUMP{u}{\rm t}}{\JUMP{u}{\rm n}} \ \text{ and }\  \tan{\psi_{\sigma}}=\frac{\sigma_t}{\sigma_n}=\frac{\kappa_{\rm t}\JUMP{u}{\rm t}}{\kappa_{\rm n}\JUMP{u}{\rm n}}.
\end{align}
Thus, the following relations  hold:
\begin{align}
|\tan{\psi_{\sigma}}|=\sqrt{\frac{\kappa_{\rm t}}{\kappa_{\rm n}}}\tan{\psi_G} \ \text{ and } \ |\tan{\psi_u}|=\sqrt{\frac{\kappa_{\rm n}}{\kappa_{\rm t}}}\tan{\psi_G}.
\end{align}
It is assumed that a crack propagates if the ERR $G$ reaches the fracture energy $G_c$, that means:
\begin{align}\label{CrackInit}
\frac{\kappa_{\rm n}\JUMPPOWER{u}{\rm n}{2}}{2}+\frac{\kappa_{\rm t}\JUMPPOWER{u}{\rm t}{2}}{2}=G_c.
\end{align}
\begin{figure}
  \centering
   \def\svgwidth{0.75\columnwidth}
   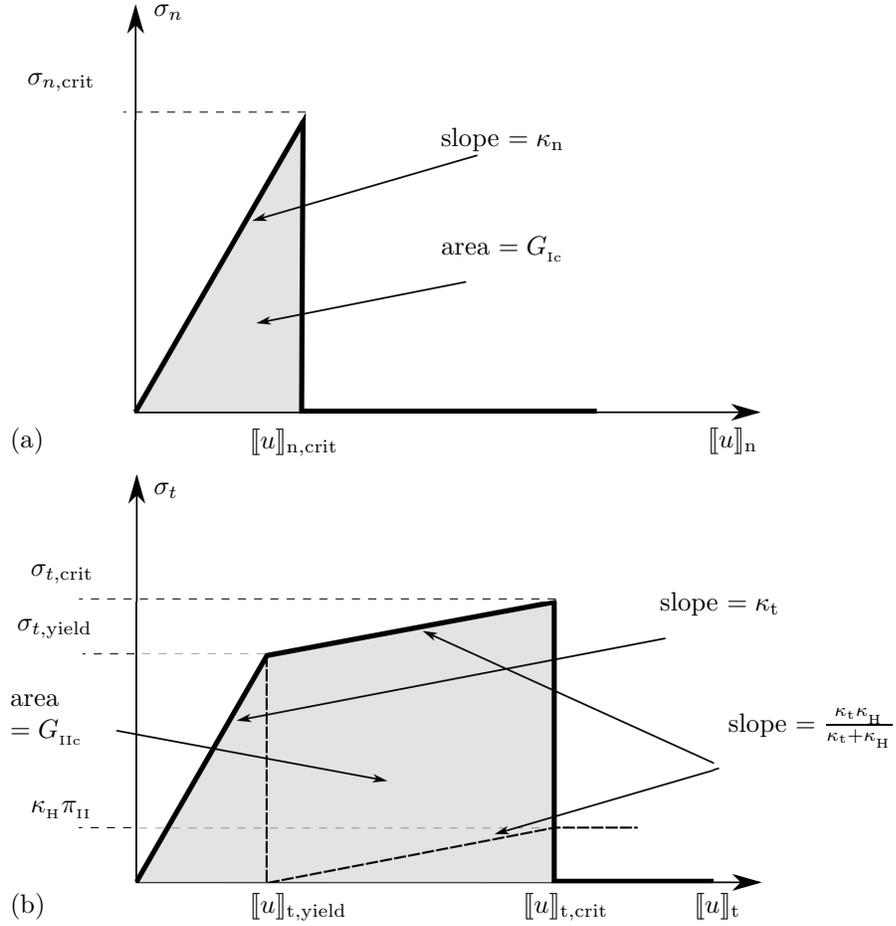
   \vspace*{2em}
   \caption{Schematic illustration of the traction-relative displacement law in the model.
   (a) pure normal (opening) mode and (b) pure tangential (shear) mode, considering $\zeta_0=1$ and $\pi_0=0$. Contribution of the delamination-gradient term is neglected, i.e. $\kappa_0{=}0$.}
\label{fig:springs}
\end{figure}

A strong dependence of $G_c$ on the fracture mode mixity  has been observed in extensive experiments \cite{BanAsk00NFCI,ERDC90FEBI,HutSuo92MMCL,LieCha92ASIF}. In accordance with other experimental observations \cite{EXPPLSL2,TveHut93IPMM}, the associated plastic zones in the adjacent bulk, near the crack tip, are larger in Mode II than in Mode I and these plastic phenomena are localized in a relatively narrow plastic zone in the bulk in the interface vicinity. In order to provide a better representation of these experimental results, a plastic tangential slip variable $\pi$
has been introduced at the interface, which allows us, firstly, to
distinguish between fracture Mode I and II in the sense that some
additional dissipated energy is associated to interface fracture in
Mode II, and secondly to simulate these narrow plastic zones.
In such a case we can model an inelastic behaviour in the tangential response of the interface, while the response in the normal direction remains linear elastic, as shown in Figure~\ref{fig:springs}.
An engineering insight into the present interface constitutive law can be summarized by the two conditions which activate the two inelastic processes included in  the formulation \cite{RoMaPa??QMMD}. The first one is the activation criterion for damage initiation  which, for the case of $\kappa_0{=}0$, reads as
\begin{align}\label{debonding-activation}
\frac12\Big(\kappa_{\rm n}\JUMPPOWER{u}{\rm n}{2}+
\kappa_{\rm t}\big(\JUMP{u}{\rm t}\!{-}\pi\big)^2\Big)= \GIC,
\end{align}
where the left hand side represents the elastic energy stored in the adhesive.
The second one concerns the evolution of $\pi$ which is triggered when $|\sigma_{\rm t}-\kappa_{_{\rm H}}\pi|$ reaches
the activation threshold $\sigma_{\mathrm{t,yield}}$, and then,
\begin{align}\label{plasticity-activation}
|\zeta\kappa_{\rm t}(\JUMP{u}{\rm t}{-}\pi)-\kappa_{_{\rm H}}\pi| =\sigma_{\mathrm{t,yield}}.
\end{align}
A more detailed analysis of the model may be found in \cite{RoMaPa??QMMD}. The model produces the desired results if
\begin{align}\label{2-sided-condition}
\frac12\sqrt{2\kappa_{\rm t}\GIC}<\sigma_{\mathrm{t,yield}}
\le\sqrt{2\kappa_{\rm t}\GIC}.
\end{align}
 The upper bound of yield stress is necessary for making possible to initiate plastic slip before the total interface damage, while the lower one is required to avoid plastic slip evolution at some point which has already been debonded.

Thus, the ERR of a mixed mode crack for the present model is defined by the:
\begin{align}\label{ERR_plast}
G=\frac{\kappa_{\rm n} \JUMPPOWER{u}{\rm n}{2}}{2}
+\frac{\kappa_{\rm t}\big(\JUMP{u}{\rm t}{-}\pi \big)^2}{2}+\sigma_{\mathrm{t,yield}}\big|\pi \big|+\frac{\kappa_{\scriptscriptstyle\textrm{H}}\pi^2}{2},
\end{align}
where it may be seen that,
referring to Eq.\,\eqref{ERR},
ERR
is here augmented by terms concerning inelastic slip $\pi$.

In the following we will try to determine the dependence of $G_c$ on the fracture mode mixity angle $\psi_{u}$ similarly as in \cite{ROMAPlzni11}. To accomplish this task, first we eliminate the plastic slip $\pi$, which in our kinematic-hardening model may be written, for $\pi>0$, as:
\begin{align}\label{PiElim}
\pi=\frac{\kappa_{\rm t}}{\kappa_{\rm t}+\kappa_{\scriptscriptstyle\textrm{H}}}\left(\JUMP{u}{\rm t}-\frac{\sigma_{t,\rm yield}}{\kappa_{\rm t}}\right).
\end{align}
Substituting Eq.\,\eqref{PiElim} into the damage initiation criterion of
Eq.\,\eqref{debonding-activation}, leads to the relation,
\begin{align}\label{DispRelation}
\frac12\left(\kappa_{\rm n}\JUMPPOWER{u}{\rm n}{2}{+}
\kappa_{\rm t}\Big(\frac{\kappa_{\scriptscriptstyle\textrm{H}}}{\kappa_{\scriptscriptstyle\textrm{H}}{+}\kappa_{\rm t}}\Big)^2\big(\JUMP{u}{\rm t}{+}\frac{\sigma_{t,\rm yield}}{\kappa_{\scriptscriptstyle\textrm{H}}}\big)^2\right){=}\GIC
\end{align}
when some interface plasticity occurs, i.e. $\JUMP{u}{\rm t}\geq\frac{\sigma_{t,\rm yield}}{\kappa_{\rm t}}$.

In a similar way as in Eq.\,\eqref{CrackInit}, which is valid if no plasticity has occurred, Eq.\,\eqref{DispRelation} defines the relation between the two components of the relative displacement at the crack tip leading to the crack growth, if some plasticity has  already appeared. This relation can be written in
a parameterized form through the use of a parametric angle $\phi$, as:
\begin{align}\label{param}
&\JUMP{u}{\rm n}=\sqrt{\frac{2\GIC}{\kappa_{\rm n}}}\cos{\phi},\nonumber \\
&\JUMP{u}{\rm t}=\sqrt{\frac{2\GIC}{\kappa_{\rm t}}}\frac{\kappa_{\rm t}+\kappa_{\scriptscriptstyle\textrm{H}}}{\kappa_{\scriptscriptstyle\textrm{H}}}\sin{\phi}-\frac{\sigma_{t,\rm yield}}{\kappa_{\scriptscriptstyle\textrm{H}}},
\end{align}
for
$\arcsin{\frac{\sigma_{t,\rm yield}}{\sqrt{2 \kappa_{\rm t} \GIC}}}\le \phi \le \frac{\pi}{2}$. Before plasticity occurs, i.e. for $0\le \phi \le \arcsin{\frac{\sigma_{t,\rm yield}}{\sqrt{2 \kappa_{\rm t} \GIC}}}$, the analogous parameterization writes as
\begin{align}\label{param_elasticpart}
&\JUMP{u}{\rm n}=\sqrt{\frac{2\GIC}{\kappa_{\rm n}}}\cos{\phi},\nonumber \\
&\JUMP{u}{\rm t}=\sqrt{\frac{2\GIC}{\kappa_{\rm t}}}\sin{\phi},
\end{align}
and angle $\phi$ coincides with the fracture mode mixity  angle $\psi_G$ defined in Eq.\,\eqref{FMMAE}. Parameterization of Eq.\,\eqref{param}, defines an ellipse whose center is at the point $(0,-\frac{\sigma_{t,\rm yield}}{\kappa_{\scriptscriptstyle\textrm{H}}})$, which   continuously
switches from the ellipse with the center at the origin of coordinates Eq.\,\eqref{param_elasticpart}, which corresponds to a state of zero plasticity.

The relation $G_c=G_c(\psi_u)$, for the case of non-zero interface plasticity, can be obtained by substitution of Eqs.\,\eqref{param}, \eqref{DispRelation} and \eqref{PiElim} into Eq.\,\eqref{ERR_plast}, leading after some algebra   to:
\begin{align}\label{GCDEP}
\GIC\left(1+\frac{\kappa_{\rm t}}{\kappa_{\scriptscriptstyle\textrm{H}}}\sin^2{\phi} \right)
-\frac{\sigma^2_{t,\rm yield}}{2\kappa_{\scriptscriptstyle\textrm{H}}}=G_c(\phi).
\end{align}
Finally, finding the relation between the angles $\phi$ and   $\psi_u$
\begin{align}
\tan{\psi_u}=\sqrt{\frac{\kappa_{\rm n}}{\kappa_{\rm t}}}
\frac{\kappa_{\rm t}{+}\kappa_{\scriptscriptstyle\textrm{H}}}{\kappa_{\scriptscriptstyle\textrm{H}}}\tan{\phi}
-\sqrt{\frac{\kappa_{\rm n}}{2\GIC}}\frac{\sigma_{t,\rm yield}}{\kappa_{\scriptscriptstyle\textrm{H}}}
\frac{1}{\cos{\phi}},
\end{align}
we obtain the desired relation $\phi=\phi(\psi_u)$ to be substituted into
Eq.\,\eqref{GCDEP}. However, an explicit relation of $G_c(\psi_u)$ is rather
cumbersome. Nevertheless, according to plots presented in \cite{ROMAPlzni11} the functional dependence of $G_c(\psi_u)$ qualitatively represents the expected behaviour in view of the previous experimental results \cite{BanAsk00NFCI,ERDC90FEBI,HutSuo92MMCL,LieCha92ASIF}.

\section{Numerical implementation} \label{sec:numimpl}
The theoretical framework, briefly presented up to this point,
provides an implementable and efficient numerical scheme.
An emerging global minimization problem, inherent in
Eq.\,\eqref{stability} may be defined by an implicit time discretization.
By discretizing the time incremental formulation in space by some
appropriate method, the problem may be casted
in a standard algebraic form.
Since the problem may be (and here is) formulated on the boundary, as all nonlinear processes considered occur exclusively on the  boundaries,
$\GC$
only,
a boundary element method seems to be a natural approach especially if
the bulk equations can efficiently be solved, which is, in particular,
the case of isotropic linear elastic materials considered in this article.
Such a formulation was developed in \cite{RoMaPa??QMMD}, using the collocation BEM but without providing a thorough description of the numerical implementation. A related symmetric Galerkin SGBEM formulation can be found in \cite{ROMAPlzni11} and a FEM implementation  in \cite{RoKrZe11DACM}. Preliminary comparison with the SGBEM formulation has shown an excellent agreement  in a few specific case studies. An advantage of the present approach with respect to a related FEM approach is that no  bulk discretization is required here, and in the analysis and optimization procedures we directly work with a relatively small number of variables associated to boundaries  in particular to $\GC$.

\subsection{Minimization problem}
Making an implicit time discretization by adopting, for simplicity,
an equidistant partition of $[0,T]$ with a fixed
time-step $\tau>0$, assuming $T/\tau\in\N$, Eq.\,\eqref{stability} leads to
a recursive minimization problem:
\begin{align}
\label{increment}
\left.
\begin{array}{ll}
\text{minimize}&\calF^k(\uC,z)=\calE(k\tau,\uC,z)+\calR(z{-}z^{k-1})\\
\text{subject to}\ &B_I \uC{\geq} 0,\quad
0{\leq}\zeta{\leq}\zeta^{k-1},
\end{array}\right\}
\end{align}
to be solved successively for $k=1,...,T/\tau$, starting from $u_0$
and $z_0$. Operator $B_I$ represents the non-pen\-etra\-tion Signorini
conditions, while the \COL{further const\-raint in \eqref{increment}}
refers to non-negativity and
irreversibility of damage parameter evolution. According to
the convention of Section~\ref{subsect-dissip},
 only $\uC$, the displacement at interfaces (or contact zones), appears in Eq.~\eqref{increment}    making clear that only this part of the displacement field is a minimizer of the problem.

We denote by $(\uC^k,z^k)$ some (generally not unique) solution to
the problem \eqref{increment}.

In order to numerically solve the emerging minimization problems
\eqref{increment}, we have utilized and test in this work several
algorithms, such as the L-BFGS-B \cite{LBFGSB} for general large scale simply
bounded problems, the GLPK routines for linear programming problems \cite{GLPK}
as well as a conjugate gradient based algorithm with constraints, see
\cite{Do09OQPA}, for solving quadratic programming problems.
Notably, the minimization problem appears to have an $L_1$-type
non-smooth term with respect to the plastic tangent slip variable $\pi$, see
Eq.\,\eqref{DissipationAdh}. In order to overcome this difficulty we take advantage of \emph{gradient projection} algorithm presented in \cite{NonSmooth} for such kind of non-smoothness.

\subsubsection{Alternate minimization algorithm}
The functional $\calF^k$ in Eq.\,\eqref{increment} is not convex
and as such   leads to a difficult minimization problem. In order to overcome this difficulty we utilize a special technique, originally proposed in \cite{BoFrMa00NERB}, called as \emph{alternate minimization algorithm} (AMA). The AMA procedure, in our case, consists in splitting the original nonconvex minimization problem to two distinct convex problems with respect to the kinematical variables ($u, \pi$) and to damage variable $\zeta$, respectively. Convergence is succeeded through an iterative procedure by alternation of this two convex problems. A flowchart of AMA may be seen in Table \ref{table:AMA}. It is worth mentioning that the individual sub-problems emerging by using such alternation consist of a nonsmooth quadratic programming problem, step (2-b), and a linear programming problem, step (2-c) of Table \ref{table:AMA}, respectively, for which we may use appropriate specialized algorithms such as those mentioned above.

\begin{table}
\begin{center}
\caption{Pseudocode of the alternate minimization algorithm}\label{table:AMA}
\framebox[1.05\width]{\minibox{
(1) \hspace{2 mm} Set $j=0$ and $\zeta^{0}=\zeta^{k-1}$ \\
(2) \hspace{2 mm} Repeat \\
\hspace{5 mm} (a) \hspace{2 mm} Set $j=j+1$\\
\hspace{5 mm} (b) \hspace{2 mm} Solve for $\uC^{j}$ and $\pi^{j}$: \\
\hspace{8 mm} minimize \hspace{2 mm} $(\uC^j,\pi^j)\mapsto\calE(t,\uC^j,z^j)+\calR(z^j{-}z^{k-1})$ \\
\hspace{8 mm} subject to \hspace{2 mm} \!\!$B_I \uC^{j}{\geq} 0$\\
\hspace{5 mm} (c) \hspace{2 mm}  Solve for $\zeta^{j}$:\\
\hspace{8 mm} minimize \hspace{2 mm} $\zeta^{j} \mapsto \calE(t,\uC^{j},z^{j})+\calR(z^{j}{-} z^{k-1})$ \\
\hspace{8 mm} subject to \hspace{2 mm} \!\!$\zeta^{k-1} \geq \zeta^{j} \geq 0$ \\
\hspace{5 mm} (d) \hspace{2 mm} If  $\parallel \zeta^{j} {-} \zeta^{j-1} \parallel < \epsilon$ exit loop \\
(3) \hspace{2 mm} Set $\uC^{k}=\uC^{j}$ and $z^{k}=z^{j}$
}}
\end{center}
\end{table}

\subsubsection{Back-tracking technique}
The above
AMA procedure
does not necessarily lead to a
globally minimizing
 solution which is, however, one of
the main
ingredient behind the
energetic-solution concept, as shown in
Section~\ref{sec:Theory}. In order to
 execute the global minimization more successfully  at particular time levels
we
use heuristic back-tracking algorithm (BTA),
devised and tested on such sort of problems
in \cite{Bene09MSMA,Bene??GONS,MiRoZe10CDEV,RoKrZe11DACM,RoMaPa??QMMD}.
\COL{The} BTA technique \COL{is} based on   checking
a two-sided energy estimate, the integral expression in Table~\ref{table:BTA},
where also some pseudo-code of BTA is given. This two-sided inequality
has been constructed by use of the energy stability condition Eq.\,\eqref{stability} and
a full deduction of it
can be found in  \cite{Miel05ERIS,MieRou09NARI,RoKrZe11DACM}. The upper and
lower energy estimates are given as time integrals of the power while
$(\overline{u}_{\rm c},\overline{z})$ and
$(\underline{u}_{\rm c},\underline{z})$
are piecewise constant interpolants in time  defined by
\begin{subequations}\label{constant}
\begin{align}
&&&\overline{u}_{\rm c}(t)=\uC^k&&\text{ for }t\in\big((k{-}1)\tau,k\tau\big],&&&&
\\
&&&\underline{u}_{\rm c}(t)=\uC^{k-1}&&
\text{ for }t\in\big[(k{-}1)\tau,k\tau\big).&&&&
\end{align}
\end{subequations}
Similar notation concerns also $\overline z$ and $\underline z$. A thorough
deduction of the boundary forms for these integrals of power,
amenable into the boundary element context, is given in
Sections~\ref{subsec:BEM} and \ref{subsec:BoundFormPotEnergy}.
\begin{table}
\begin{center}
\caption{Pseudocode
of energy-based backtracking algorithm}\label{table:BTA}
\framebox[1.01\width]{\minibox{
(1) \hspace{1 mm} Set $k=1$ and $\zeta^{0}=\zeta_0$ \\
(2) \hspace{1 mm} Repeat \\
\hspace{2 mm} (a) \hspace{1 mm} Determine $\zeta^{k}$ using the alternating minimization\\
\hspace{2 mm}  algorithm for time $t_k$ and the initial value $\zeta^{0}$\\
\hspace{2 mm} (b) \hspace{1 mm} Set $\zeta^{0}=\zeta^{k}$\\
\hspace{2 mm} (c) \hspace{1 mm} If the two-sided energy estimate holds:\\
\hspace{6 mm} $\int_{(k-1)\tau}^{k\tau}\calE_t'(t,
\COL{\overline{u}_{\rm c},\overline{z}})\,\d t \leq \calE(k\tau,\uC^{k},z^{k})+\calR(z^{k} {-} z^{k-1})$ \\
\hspace{11.5mm}$-\calE((k{-}1)\tau,\uC^{k-1},z^{k-1})\le
\int_{(k-1)\tau}^{k\tau}\calE_t'(t,\COL{\underline{u}_{\rm c},\underline{z}})\,\d t$\\
\hspace{5 mm} set $k=k+1$\\
\hspace{2 mm} (d) \hspace{1 mm} Else set $k=k-1$ \\
\hspace{2 mm} (e) \hspace{1 mm} Until $k > T/\tau$
}}
\end{center}
\end{table}
Although there is no proof that BTA converges to the global minimum,
definitely it leads to solutions of lower energy than those obtained if we
\COL{used AMA only}.

\subsection{Boundary element method}\label{subsec:BEM}
\COL{The} boundary element method is closely related to the map between the
prescribed boundary conditions in displacements or tractions and the unknown
boundary  displacements or tractions. In pure Dirichlet and Neumann
boundary\COL{-}value problems (BVPs), these maps are called
\emph{Steklov-Poincar\'e} and  \emph{Poincar\'e-Steklov maps} \cite{KhoWitt2004,SaSch10}, respectively,
and BEM can be considered as an approach to discretize these maps.
In the present computational procedure, the role of the BEM analysis, applied to each subdomain $\Omega_i$ separately (which, in fact, makes this problem very suitable for parallel computers), is  to solve the corresponding BVPs on each $\Omega_i$.
For this goal, we  numerically solve  the
\emph{Somigliana displacement identity} \cite{PaCa1997,SaSch10}
\begin{align}\label{eq_Somigliana}
c^i_{ml}(\xi)u^i_m(\xi){+}\int_{\Gamma_i} \!\!\!\!\!\!\!\!\!
-\ u^i_m(x)T^i_{ml}(x,\xi)\,\d S_x\nonumber \\
{=}\int_{\Gamma_i} p^i_m(x)U^i_{ml}(x,\xi)\,\d S_x,
\end{align}
where $\xi \in \Gamma_i=\partial\Omega_i$ \COL{and} $u^i_m(x)$ and  $p^i_m(x)$ denote
the $m$-component of the displacement and traction vector, respectively. The superscript $i$ used in this section refers to domain $\Omega_i$ in difference to previous and next sections where it denotes time step.  The weakly singular integral kernel $U^i_{ml}(x,\xi)$, two-point tensor field, given by the Kelvin fundamental solution  (free-space Green's function) represents the displacement at $x$ in the $m$-direction originated by a unit point force  at $\xi$ in  the $l$-direction  in the unbounded elastic medium whose material properties coincide with those of $\Omega_i$. The strongly singular integral kernel $T^i_{ml}(x,\xi)$, two-point tensor field, represents the corresponding tractions at $x$ in the $m$-direction. The co\-ef\-fi\-cient-tensor $c^i_{ml}(\xi)$ of the free-term is a function of the local geometry of the boundary $\Gamma_i$ at $\xi$, and may be evaluated by a closed analytical formula for  isotropic elastic solids \cite{Mantic1993CM}. The symbol $\int \!\!\!\!\!\!\,-$
in Eq.\,\eqref{eq_Somigliana}
stands for the Cauchy principal value of an integral.

Consider a discretization of the boundary $\Gamma_i$ by a boundary element mesh, which is also used to define a suitable  discretization of boundary displacements $u^i(x)$ and tractions $p^i(x)$ by interpolations of their nodal values. By imposing (collocating)  the Somigliana identity \eqref{eq_Somigliana} at all boundary nodes (called collocation points)  we set the BEM
system of linear equations for $\Gamma_i$. The solution of this system defines the unknown nodal values of displacements and tractions along $\Gamma_i$ representing a part of arrays of all nodal values denoted as $\mathbf{u}^i$ and $\mathbf{p}^i$, respectively.
The arrays $\mathbf{u}^i$ and $\mathbf{p}^i$ also include the known nodal values of   displacements and tractions along $\Gamma_i$ given by the prescribed boundary conditions.
The BEM system obtained from Eq.\,\eqref{eq_Somigliana} is usually written as $\mathbf{H}^i\mathbf{u}^i=\mathbf{G}^i\mathbf{p}^i$~\cite{PaCa1997}.
In our computer implementation of BEM, we employ straight
elements with continuous and piecewise linear \COL{interpolation} for displacements
and  possibly discontinuous piecewise linear interpolation for tractions.

Then, to compute an approximation of the  elastic energy,  $\calE_{\rm el}$ from Eq.\,\eqref{Eell}, stored in each bulk $\Omega_i$, by using the obtained approximations of boundary displacements $u_m^i$ and  of the corresponding boundary tractions $p_m^i$ along $\Gamma_i$, we utilize the following general relation \cite{Hartmann}, neglecting body forces:
\begin{align}\label{eq_ENERGYIDENTITY}
& {\calE}_{\Omega_i} (t,u^i):=\!\left\{\begin{array}{ll}
  \displaystyle{\
\!\frac{1}{2}\int_{\Gamma_i}\!\!\!\! u^i p^i(u^i) }
\displaystyle{\,\d S}
&\text{ if }u^i{=}{\uDi}(t) \text{ on }{\GDi},\\[.3em]
   \infty&\text{ elsewhere,}
\end{array}\right.\!\!\!\!
\end{align}
while the corresponding total potential energy  is
\begin{equation}\label{eq_TPE}
\Pi_{\Omega_i} (t,u^i)
={\calE}_{\Omega_i} (t,u^i)-\int_{{\GNi}} u^i {\pNi}\,\d S.
\end{equation}
Notice that both Eq.\,\eqref{eq_ENERGYIDENTITY} and  \eqref{eq_TPE} provide \emph{pure boundary expressions of energy}. If  $\calE_{\rm el}$ in Eq.\,\eqref{StoredEnergyDef} is replaced by the sum of potential energies $\Pi_{\Omega_i}$ for all subdomains, the functional $\calE$  defined by Eq.\,\eqref{StoredEnergyDef}  will represent the total potential energy of the whole problem.

We will also need to compute integrals of time derivatives of energy, appearing in Table~\ref{table:BTA}, were it is presupposed that displacements on the contact boundary part $\uC$ are defined in time steps $k$ and $k{-}1$. In order to do such calculations we need to separate the problem into three different sub-problems, in each of them either the contact or prescribed Dirichlet or Neumann data are defined on the boundary. This separation to sub-problems may also serve to express and solve the minimization problem, considering the integral on   $\GC$   only.

\subsection{Boundary forms of the total potential energy for a single domain}\label{subsec:BoundFormPotEnergy}
%
\begin{figure*}
   \centering
   \def\svgwidth{1.0\textwidth}
   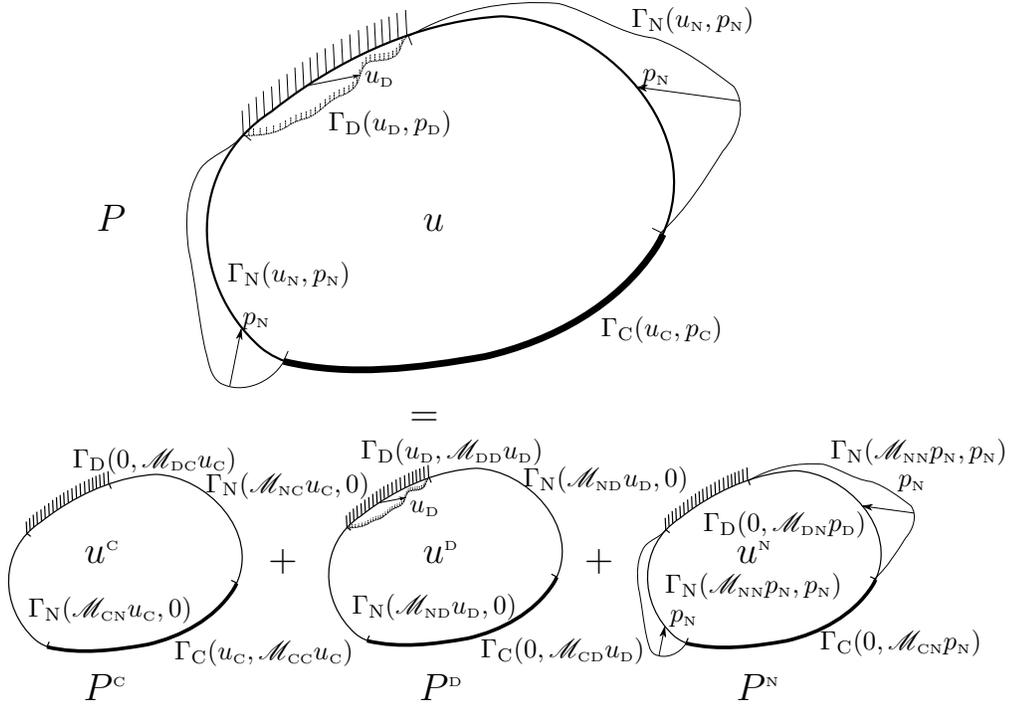
   \hspace*{.3em}
   \caption{Solution of a mixed BVP  $\SO$ for a single elastic domain given as a superposition
of the solutions of the three sub-problems $\SC,\SD$ and $\SN$.}
\label{fig:states}
\end{figure*}
Consider the BVP for a sub-domain $\Omega_i$. In this section we will omit index $i$, for the sake of simplicity. Let $u_\eta$ and $p_\eta$, respectively,  denote the displacement   and traction solutions of this BVP restricted to $\Gamma_\eta$, $\eta=\text{\rm C}$, $\text{\rm D}$ and $\text{\rm N}$,
 e.g. $\uC=u|_{\GC}$ and $\pD=p|_{\GD}$.
We assume here a mixed-type operator $\calM$
which formally assigns $(\pC,\pD,\uN)$
 to the known boundary data $(\uC,\uD,\pN)$
of the original BVP  $\SO$ shown in Fig.~\ref{fig:states}, and may be
expressed using the following block structure as:
\begin{align}\label{Moperator}
\Bigg(\!\begin{array}{l}\pC\\\pD\\\uN\end{array}\!\Bigg)=
\Bigg(\!\begin{array}{ccc}
\calM_{\text{\tiny\rm CC}}&\calM_{\text{\tiny\rm CD}}&\calM_{\text{\tiny\rm CN}}\\
\calM_{\text{\tiny\rm DC}}&\calM_{\text{\tiny\rm DD}}&\calM_{\text{\tiny\rm DN}}\\
\calM_{\text{\tiny\rm NC}}&\calM_{\text{\tiny\rm ND}}&\calM_{\text{\tiny\rm NN}}
\end{array}\!\Bigg)
\Bigg(\!\begin{array}{l}\uC\\\uD\\\pN\end{array}
\!\Bigg).
\end{align}
The columns of the aforementioned block operator $\calM$ are associated to the   sub-problems  $P^\eta$ defined in Fig.~\ref{fig:states}. The displacement solution of a subproblem $P^\eta$ is denoted as $u^\eta$.
 From the principle of superposition the  displacement solution of $\SO$ may be reconstructed by the sum:
\begin{align}\label{eq:totDisp}
u=\uCS+\uDS+\uNS \,,
\end{align}

The total potential energy for the mixed type BVP $\SO$ can be written in an expanded form as,
\begin{align}\label{TPE_SO}
\Pi (t,u){=}\frac{1}{2}\!\int_{\GC}\!\!\!\! \uC \pC \d S{+}\frac{1}{2}\!\int_{\GD}\!\!\!\! \uD \pD \d S {-}\frac{1}{2}\!\int_{\GN}\!\!\!\! \uN \pN \d S.
\end{align}
By substituting the unknown data for the problem $\SO$ from Eq.\,\eqref{Moperator}
the total potential energy   writes  as
\begin{align}\label{extendTPE_SO}
&\Pi (t,\uC)=\frac{1}{2}\bigg(\int_{\GC}\!\!\!\! \uC \calM_{\text{\tiny\rm CC}}\uC \,\d S
\nonumber \\
&\hspace{7.6em}
{+}\int_{\GC}\!\!\!\! \uC \calM_{\text{\tiny\rm CD}}\uD \,\d S{+}\int_{\GC}\!\!\!\! \uC \calM_{\text{\tiny\rm CN}}\pN \,\d S
\nonumber \\
&{+}\int_{\GD}\!\!\!\! \uD \calM_{\text{\tiny\rm DC}}\uC \,\d S{+}\int_{\GD}\!\!\!\! \uD \calM_{\text{\tiny\rm DD}}\uD \,\d S{+}\int_{\GD}\!\!\!\! \uD \calM_{\text{\tiny\rm DN}}\pN \,\d S
\nonumber \\
&{-}\int_{\GN}\!\!\!\! \pN \calM_{\text{\tiny\rm NC}}\uC \,\d S{-}\int_{\GN}\!\!\!\! \pN \calM_{\text{\tiny\rm ND}}\uD \,\d S{-}\int_{\GN}\!\!\!\! \pN \calM_{\text{\tiny\rm NN}}\pN \,\d S\bigg).
\end{align}

From Eq.\,\eqref{extendTPE_SO} it is clear, that since $\uD(t)$ and $\pN(t)$ are  known, the total potential energy is in fact  a function of   the contact displacement  $\uC$ only,  in addition to be a function of time $t$. We further modify Eq.\,\eqref{extendTPE_SO}, in order to hold the unknown variables  on $\GC$ only, by utilizing   \emph{the second Betti reciprocity relation} between the elastic solutions of  $\SC$ and   $\SN$,
\begin{align}\label{REC_C_N}
\int_{\GN}\!\!\!\! \pN \calM_{\text{\tiny\rm NC}}\uC \,\d S =- \int_{\GC}\!\!\!\! \uC \calM_{\text{\tiny\rm CN}}\pN \,\d S
\end{align}
as well as between the solutions of $\SC$ and $\SD$,
\begin{align}\label{REC_C_D}
\int_{\GD}\!\!\!\! \uD \calM_{\text{\tiny\rm DC}}\uC \,\d S = \int_{\GC}\!\!\!\! \uC \calM_{\text{\tiny\rm CD}}\uD \,\d S.
\end{align}
Then, by substituting Eqs.\,\eqref{REC_C_N} and \eqref{REC_C_D} into Eq.\,\eqref{extendTPE_SO},
\begin{align}\label{extendTPE_SO_modified}
&\Pi (t,\uC)=
\int_{\GC}\!\!\!\! \uC \Big(\frac{1}{2}\calM_{\text{\tiny\rm CC}}\uC{+}\calM_{\text{\tiny\rm CD}}\uD{+}\calM_{\text{\tiny\rm CN}}\pN \Big)  \,\d S
\nonumber \\
&\hspace{10mm}{+}\frac{1}{2}\int_{\GD}\!\!\!\! \uD \calM_{\text{\tiny\rm DD}} \uD \,\d S{+}\frac{1}{2}\int_{\GD}\!\!\!\! \uD \calM_{\text{\tiny\rm DN}}\pN \,\d S
\nonumber \\
&\hspace{10mm}{-}\frac{1}{2}\int_{\GN}\!\!\!\! \pN \calM_{\text{\tiny\rm ND}}\uD \,\d S{-}\frac{1}{2}\int_{\GN}\!\!\!\! \pN \calM_{\text{\tiny\rm NN}}\pN \,\d S.
\end{align}

  The partial time derivative of the total potential energy expression in Eq.\,\eqref{extendTPE_SO_modified}  writes as
\begin{align}\label{extendTPE_SO_dt}
&\frac{\partial \Pi}{\partial t} (t,\uC)=
\int_{\GC}\!\!\!\! \uC \big(\calM_{\text{\tiny\rm CD}}\duD{+}\calM_{\text{\tiny\rm CN}}\dpN \big)  \,\d S
\nonumber \\
&\hspace{10mm}{+}\frac{1}{2}\int_{\GD}\!\!\!\!
\DT{\overline{\uD \calM_{\text{\tiny\rm DD}}\uD}} \,\d S
{+}\frac{1}{2}\int_{\GD}\!\!\!\!
\DT{\overline{\uD \calM_{\text{\tiny\rm DN}}\pN}} \,\d S
\nonumber \\
&\hspace{10mm}{-}\frac{1}{2}\int_{\GN}\!\!\!\! \DT{\overline{\pN \calM_{\text{\tiny\rm ND}}\uD}} \,\d S{-}\frac{1}{2}\int_{\GN}\!\!\!\! \DT{\overline{\pN \calM_{\text{\tiny\rm NN}}\pN}} \,\d S
\end{align}
where the bar with dot denotes the time derivative of the expression below the bar.
The integral corresponding to that on the left-hand side in the two-sided inequality in Table~\ref{table:BTA}, can be evaluated using Eq.\,\eqref{extendTPE_SO_dt} as,
\begin{align}\label{LowerTSI}
&\int_{(k-1)\tau}^{k\tau}\frac{\partial \Pi}{\partial t} (t,\uC^k)=
\int_{\GC}\!\!\!\! \uC^{k} \big(\calM_{\text{\tiny\rm CD}}\uD^{k}{+}\calM_{\text{\tiny\rm CN}}\pN^{k} \big)  \,\d S
\nonumber \\
&\hspace{10mm} {-}\int_{\GC}\!\!\!\! \uC^{k} \big(\calM_{\text{\tiny\rm CD}}\uD^{k-1}{+}\calM_{\text{\tiny\rm CN}}\pN^{k-1} \big)  \,\d S
\nonumber \\
&\hspace{10mm}{+}\frac{1}{2}\int_{\GD}\!\!\!\! \uD^{k} \calM_{\text{\tiny\rm DD}}\uD^{k} \,\d S{-}\frac{1}{2}\int_{\GD}\!\!\!\! \uD^{k-1} \calM_{\text{\tiny\rm DD}}\uD^{k-1} \,\d S
\nonumber \\
&\hspace{10mm}{+}\frac{1}{2}\int_{\GD}\!\!\!\! \uD^{k} \calM_{\text{\tiny\rm DN}}\pN^{k} \,\d S{-}\frac{1}{2}\int_{\GD}\!\!\!\! \uD^{k-1} \calM_{\text{\tiny\rm DN}}\pN^{k-1} \,\d S
\nonumber \\
&\hspace{10mm}{-}\frac{1}{2}\int_{\GN}\!\!\!\! \pN^{k} \calM_{\text{\tiny\rm ND}}\uD^{k} \,\d S{+}\frac{1}{2}\int_{\GN}\!\!\!\! \pN^{k-1} \calM_{\text{\tiny\rm ND}}\uD^{k-1} \,\d S
\nonumber \\
&\hspace{10mm}{-}\frac{1}{2}\int_{\GN}\!\!\!\! \pN^{k} \calM_{\text{\tiny\rm NN}}\pN^{k} \,\d S{+}\frac{1}{2}\int_{\GN}\!\!\!\! \pN^{k-1} \calM_{\text{\tiny\rm NN}}\pN^{k-1} \,\d S,
\end{align}
and similarly for the integral on the right-hand side of the two-sided inequality in Table~\ref{table:BTA},
\begin{align}\label{UpperTSI}
&\int_{(k-1)\tau}^{k\tau}\frac{\partial \Pi}{\partial t} (t,\uC^{k-1})=
\int_{\GC}\!\!\!\! \uC^{k-1} \big(\calM_{\text{\tiny\rm CD}}\uD^{k}{+}\calM_{\text{\tiny\rm CN}}\pN^{k} \big)  \,\d S
\nonumber \\
&\hspace{10mm}{-}\int_{\GC}\!\!\!\! \uC^{k-1} \big(\calM_{\text{\tiny\rm CD}}\uD^{k-1}{+}\calM_{\text{\tiny\rm CN}}\pN^{k-1} \big)  \,\d S
\nonumber \\
&\hspace{10mm}{+}\frac{1}{2}\int_{\GD}\!\!\!\! \uD^{k} \calM_{\text{\tiny\rm DD}}\uD^{k} \,\d S{-}\frac{1}{2}\int_{\GD}\!\!\!\! \uD^{k-1} \calM_{\text{\tiny\rm DD}}\uD^{k-1} \,\d S
\nonumber \\
&\hspace{10mm}{+}\frac{1}{2}\int_{\GD}\!\!\!\! \uD^{k} \calM_{\text{\tiny\rm DN}}\pN^{k} \,\d S{-}\frac{1}{2}\int_{\GD}\!\!\!\! \uD^{k-1} \calM_{\text{\tiny\rm DN}}\pN^{k-1} \,\d S
\nonumber \\
&\hspace{10mm}{-}\frac{1}{2}\int_{\GN}\!\!\!\! \pN^{k} \calM_{\text{\tiny\rm ND}}\uD^{k} \,\d S{+}\frac{1}{2}\int_{\GN}\!\!\!\! \pN^{k-1} \calM_{\text{\tiny\rm ND}}\uD^{k-1} \,\d S
\nonumber \\
&\hspace{10mm}{-}\frac{1}{2}\int_{\GN}\!\!\!\! \pN^{k} \calM_{\text{\tiny\rm NN}}\pN^{k} \,\d S{+}\frac{1}{2}\int_{\GN}\!\!\!\! \pN^{k-1} \calM_{\text{\tiny\rm NN}}\pN^{k-1} \,\d S.
\end{align}
For the case where  homogeneous boundary conditions are prescribed on $\GN$, i.e. $\pN=0$, the above equations are simplified further. In such a case the total potential energy will coincide with the \emph{elastic strain energy} and Eq.\,\eqref{extendTPE_SO_modified} takes the form,
\begin{align}\label{extendTPE_SO_modified_mod}
\Pi (t,\uC)&=
\int_{\GC} \uC \Big(\frac{1}{2}\calM_{\text{\tiny\rm CC}}\uC{+}\calM_{\text{\tiny\rm CD}}\uD \Big)  \,\d S
\nonumber \\
&+\frac{1}{2}\int_{\GD} \uD \calM_{\text{\tiny\rm DD}} \uD \,\d S,
\end{align}
In this case, the time derivative, given by Eq.\,\eqref{extendTPE_SO_dt}, is written as,
\begin{align}\label{extendTPE_SO_dt_hom}
\frac{\partial \Pi}{\partial t} (t,\uC){=}
\int_{\GC}\!\!\uC \big(\calM_{\text{\tiny\rm CD}}\duD \big)  \,\d S{+}\frac{1}{2}\int_{\GD}\!\!\DT{\overline{\uD \calM_{\text{\tiny\rm DD}}\uD}} \,\d S.
\end{align}

Furthermore, in this case $(\pN=0)$, the lower and upper energy estimates in the two-sided inequality  are further simplified as
\begin{align}\label{LowerTSI_hom}
&\int_{(k-1)\tau}^{k\tau}\!\!\frac{\partial \Pi}{\partial t} (t,u^{k})=
\int_{\GC}\! \uC^{k} \big(\calM_{\text{\tiny\rm CD}}\uD^{k}{-}\calM_{\text{\tiny\rm CD}}\uD^{k-1} \big)  \,\d S
\nonumber \\
&\qquad{+}\frac{1}{2}\int_{\GD}\! (\uD^{k}{-}\uD^{k-1}) \big(\calM_{\text{\tiny\rm DD}}\uD^{k}{+}\calM_{\text{\tiny\rm DD}}\uD^{k-1} \big)  \,\d S,
\end{align}
and
\begin{align}\label{UpperTSI_hom}
&\int_{(k-1)\tau}^{k\tau}\!\!\frac{\partial \Pi}{\partial t} (t,u^{k-1})=
\int_{\GC}\! u^{k-1} \big(\calM_{\text{\tiny\rm CD}}\uD^{k}{-}\calM_{\text{\tiny\rm CD}}\uD^{k-1} \big)  \,\d S
\nonumber \\
&\qquad{+}\frac{1}{2}\int_{\GD}\! (\uD^{k}{-}\uD^{k-1}) \big(\calM_{\text{\tiny\rm DD}}\uD^{k}{+}\calM_{\text{\tiny\rm DD}}\uD^{k-1} \big)  \,\d S,
\end{align}
respectively.

The above formulation assumes the solution of each of the three sub-problems in all time steps, while the total response results by their superposition. It is worth mentioning that for proportional external loading, which is  separable with respect to spatial coordinates  and time, that means
\begin{subequations}\label{separable}
\begin{align}
&\uD (t,x)= \phi(t)\,\uD(x)\ \ \ \text{ for } x
\COL{\in}\GD,\label{tfu}    \\
&\pN (t,x)= \psi(t)\,\pN(x)\ \ \ \text{ for } x
\COL{\in}\GN,\label{tfp}
\end{align}
\end{subequations}
the $\SD$ and $\SN$ problems need to be solved just once for the first time step, while for the subsequent steps the solutions will be generated  by using functions $\phi$ and $\psi$ from Eq.\,\eqref{separable} as appropriate multipliers.

Summarizing,  expressions  where the minimizer of the total potential energy is the displacement field defined  at the adhesive  contact  boundary part $\GC$  have been established in this section.   Moreover,    appropriate formulas for computing the lower and upper energy estimates shown in Table~\ref{table:BTA} have been given.

\subsection{Interface elements}\label{subsec:InterElements}
The interconnection of the subdomains as well as the consideration of Signorini kinematical conditions is at\-tained by intermediate elements referred to as interface elements. A local reference system is associated to each interface element defining a normal and a tangential component of relative displacements. In the case of contact problems of two deformable bodies where only small changes in the geometry are assumed and  conforming meshes of   elastic domains are considered along the interface, it is possible to incorporate the contact constraints on a purely nodal basis. For a general case of nodes   arbitrarily distributed along the possible contact interface between two bodies, which can occur e.g. when automatic meshing is used for two different bodies, further considerations must be taken into account about the definition of Signorini contact conditions, this case not being considered here.
The mechanical properties of springs  distributed continuously at the interface, are given by their normal and tangential stiffnesses $\kappa_{\rm n}$ and $\kappa_{\rm t}$, respectively, and additionally in the tangent direction also  by the so-called plastic modulus $\kappa_{\scriptscriptstyle\textrm{H}}$ and the factor of influence of damage $\kappa_{\rm 0}$. The shape functions used to approximate the distribution of variables at interface elements are linear and continuous for the displacements, while for the inelastic variables, $\zeta$ and $\pi$, might be constant or alternatively, continuous or discontinuous linear. In addition to the continuous distribution of springs, the interfaces and interface elements may be equipped by a ``dissipative mechanism'' whose properties are the mode\COL{-}I fracture energy $\GIC$ and the critical stress $\sigma_{\mathrm{t,yield}}$ used in \eqref{DissipationAdh}.

\section{Numerical examples} \label{sec:numexmp}
\begin{figure}
   \centering
   \def\svgwidth{.8\columnwidth}
   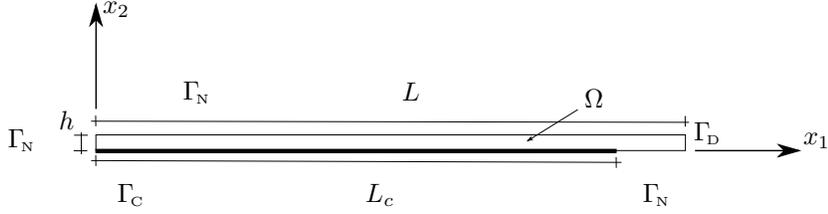
   \vspace*{2.5em}
   \caption{Problem geometry and boundary conditions.}
\label{fig_m1}
\end{figure}

The above introduced formulation has been implemented  in a  two-dimensional BEM code \cite{OBEMP} using continuous piecewise linear boundary elements \cite{PaCa1997}, and also supplied with all the necessary \emph{modules} for the EC-BEM, where the acronym EC-BEM refers to the Energetic approach for the solution of adhesive Contact problems by BEM.
The geometry of the problem solved is shown in Figure~\ref{fig_m1}. With reference to Figure~\ref{fig:geom}, only one subdomain (i.e.~$N=1$) is used to model in a simple way an experimental test motivated by the  pull-push shear test used in engineering practice \cite{CoCa11ES}. Thus, the debonding  occurs between the domain and the rigid foundation interface $\GC$.

The length and height of the rectangular domain $\Omega$, respectively, are $L=250$ mm and $h=12.5$ mm. The length of the initially glued part $\GC$ placed at the bottom side of $\Omega$ is $L_c=0.9 L=225$ mm. The \COL{isotropic}
elastic material of the bulk is
aluminium with the
\COL{Young}
modulus $E=70$ GPa and the Poisson ratio $\nu=0.35$. Elastic plain strain state is considered. The material adhesive layer is epoxy resin, with elastic properties $E_a=2.4$ GPa and $\nu_a=0.33$. Assuming the thickness of the adhesive layer $h_a=0.2$ mm, and following \cite{TMGP11BEMA}, the corresponding stiffness parameters are represented by the normal stiffness $\kappa_{\rm n}=\frac{E_a(1-\nu_a)}{h_a(1+\nu_a)(1-2\nu_a)}=$18 GPa/mm and the tangential stiffness $\kappa_{\rm t}=\frac{\kappa_{\rm n}(1-2\nu_a)}{2(1-\nu_e)}=\kappa_{\rm n}/4$. The parameters for the dissipation mechanisms are the mode\COL{-}I fracture energy $\GIC=0.01$ J/mm$^2$ as well as the yield slip stress $\sigma_{t,\rm yield}=$168 MPa. Then, $\sigma_{n,{\rm crit}}=\sqrt{2\kappa_{\rm n}\GIC}=$600MPa and $\sigma_{t,{\rm crit}}=\sqrt{2\kappa_{\rm t}\GIC}=$300MPa.
Finally, the hardening slope for plastic slip is   $\kappa_{\scriptscriptstyle\textrm{H}} =\kappa_{\rm t}/9$.

\subsection{Numerical experimentation}
A few sample problem cases are solved   in order to illustrate the capabilities of the numerical procedure developed. In Section~\ref{par:NMONOTLOAD} we  experiment  with a non-monotonic Dirichlet loading, for a variety of combination of dissipation properties. In the next Section~\ref{par:TRACTIOMLOAD} we present   results for a monotonic Neumann loading on a modified geometry of the problem in order to avoid the absence of a Dirichlet boundary part especially after the total delamination of the interface. Both   examples allow us to illustrate the behaviour of the numerical solution of the present energetic formulation for delamination problems, and do not aim to analyze the problem solutions in a thorough manner.
 In all the numerical computations, linear continuous elements have been used for the interface displacement and plastic slip variables,   also referred to as kinematical variables ($u,\pi$), while constant discontinuous elements have been assumed for the damage variable $\zeta$.

\subsubsection{Non-monotonic loading}\label{par:NMONOTLOAD}
A hard-device loading is assumed by prescribing horizontal and vertical displacements, $u_1(t,x)=\sin(7\,t) w_1(x)$, where $t \in [0,1]$ while $w_1$=1 mm, and $u_2=0.6 u_1$, respectively, at the  right-hand side  of the rectangle $\Gamma$, defining the Dirichlet boundary $\GD$. In accordance with Eq.\,\eqref{tfu}, $\phi(t){=}\sin(7\,t)$. All the other boundary parts are considered to be traction free, defining the Neumann boundary $\GN$, except for the contact surface $\GC$.
The boundary $\Gamma$ is discretized by 64 elements using a uniform boundary element mesh along each side,  27 elements being used for $\GC$.
\begin{figure}
\begin{center}
\psfrag{xxxxxxxxxxxxxxxxxxxx} {\scriptsize{Signorini contact}}
\psfrag{xxxxxxxxxxxxxxxxxxxy} {\scriptsize{Plasticity}}
\psfrag{resultant_force_x___________________________(m)} {\!\!\!\!\!\!\scriptsize{Resultant horizontal force (KN/m)}}
\psfrag{displacement_(m)} {\!\!\!\!\!\!\!\!\!\!\!\!\!\!\!\!\!\!\!\!\!\!\!\!\scriptsize{Horizontal displacement on $\GD$ (mm)}}
{\includegraphics[width=.7\textwidth, height=.50\textwidth]{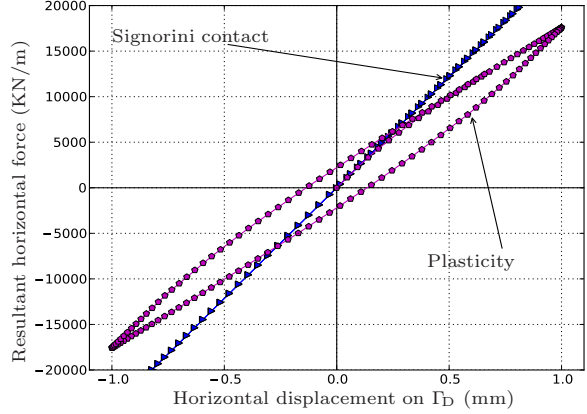}}
\caption{The horizontal resultant force  versus the horizontal displacement on $\GD$: (a)   pure Signorini contact and (b) interface plasticity included.}
\label{fig:ex1_1}
\end{center}
\end{figure}
\begin{figure}
\begin{center}
\psfrag{xxxxxxxxxxxxxxxxxxxy} {\scriptsize{Signorini contact}}
\psfrag{xxxxxxxxxxxxxxxxxxxx} {\scriptsize{Plasticity}}
\psfrag{resultant_force_x___________________________(m)} {\,\,\,\scriptsize{Shear stress (MPa)}}
\psfrag{displacement_(m)} {\!\!\!\!\!\!\!\!\!\!\!\!\!\!\!\!\!\!\!\!\!\!\!\!\scriptsize{Tangential relative displacement (mm)}}
{\includegraphics[width=.7\textwidth, height=.5\textwidth]{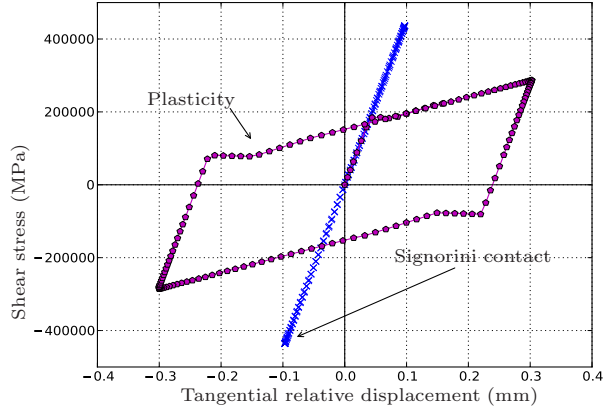}}
\caption{Stress-relative displacement behaviour computed at the interface point $x_1=208.33$ mm of $\GC$:   (a)   pure Signorini contact and (b) interface plasticity included.}
\label{fig:ex1_2}
\end{center}
\end{figure}
Four combinations of properties of the dissipative mechanism of the adhesive are considered:
\begin{enumerate}[(a)]
	\item Absence of any dissipation, leading to a pure elastic Signorini contact problem,
	\item Interface plasticity is considered, the damage variables $\zeta$ being excluded from the minimization procedure,
	\item Interface damage  is considered,  the plastic slip variables $\pi$ being excluded from the minimization procedure,
	\item Both interface damage and plasticity are considered.
\end{enumerate}
Cases (a) and (b) are mainly included for the comparison purposes and also in order to analyse   an inelastic response due to interface plasticity. The horizontal resultant force with respect to the  displacement on $\GD$ is plotted  in Fig.~\ref{fig:ex1_1}. For the case of an inelastic response due to a non-monotonic loading, a hysteresis cycle  appears as  expected. Furthermore, for these two cases also   the shear stresses with respect to the relative tangential displacements at an interface point  $x_1=208.33$ mm  are plotted in Fig.~\ref{fig:ex1_2}, where a typical hysteretic behaviour for the kinematic type hardening  plasticity is successfully computed
 in   case (b).
\begin{figure}
\begin{center}
\psfrag{xxxxxxxxxxxxxxxxxxxx} {\scriptsize{Damage}}
\psfrag{xxxxxxxxxxxxxxxxxxxy} {\scriptsize{Damage and plasticity}}
\psfrag{resultant_force_x___________________________(m)} {\!\!\!\!\!\!\scriptsize{Resultant horizontal force (KN/m)}}
\psfrag{displacement_(m)} {\!\!\!\!\!\!\!\!\!\!\!\!\!\!\!\!\!\!\!\!\!\!\!\!\scriptsize{Horizontal displacement on $\GD$ (mm)}}
{\includegraphics[width=.7\textwidth, height=.50\textwidth]{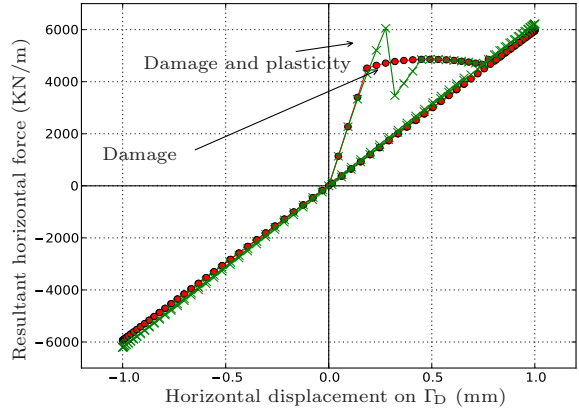}}
\caption{The horizontal resultant force   versus  the horizontal displacement on $\GD$: (c) interface damage and (d)  both interface damage and plasticity.}
\label{fig:ex1_3}
\end{center}
\end{figure}

\begin{figure}
\begin{center}
\psfrag{xxxxxxxxxxxxxxxxxxxx} {\scriptsize{Damage}}
\psfrag{xxxxxxxxxxxxxxxxxxxy} {\scriptsize{Damage and plasticity}}
\psfrag{resultant_force_y___________________________(m)} {\!\!\!\!\!\!\scriptsize{Resultant vertical force (KN/m)}}
\psfrag{displacement_(m)} {\!\!\!\!\!\!\!\!\!\!\!\!\!\!\!\!\!\!\!\!\!\!\!\!\scriptsize{Vertical displacement on $\GD$ (mm)}}
{\includegraphics[width=.7\textwidth, height=.50\textwidth]{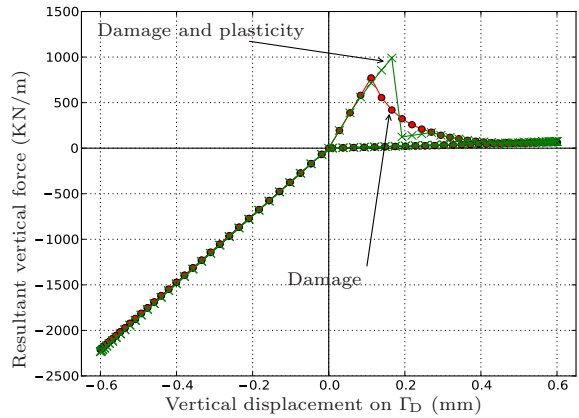}}
\caption{The vertical resultant force versus the vertical displacement on $\GD$: (c) interface damage and (d)  both  interface  damage and plasticity.}
\label{fig:ex1_4}
\end{center}
\end{figure}

More complicated behaviour is obtained for     cases (c) and (d) where   interface damage is included. This may be observed in Figs.~\ref{fig:ex1_3} and \ref{fig:ex1_4} where the horizontal and vertical resultant forces with the respective displacements at $\GD$ are depicted. For these cases upon the first uploading a damage initially appears, for   case (c) by breaking one element that corresponds to a crack opening of $L_{\rm crack}=0.037L_c=8.33$ mm, while for   case (d) in a following time step, by breaking six elements simultaneously which corresponds to a crack opening of $L_{\rm crack}=0.22L_c=50.0$ mm. For the same time step, that crack initiates in   case (d) with a crack opening of $L_{\rm crack}=50.0$ mm, for   case (c) after some progressive damage propagation,   a crack of the same length exists. This behaviour is the expected one, since because of plasticity appearance in  case (d), damage delayed to appear in comparison with case c) where  it is assumed that an energy may be released only due to damage. Then, after change of the direction of loading no further  damage  appears, while plasticity still evolves on the remaining glued part of $\GC$ upon the respective uploading periods.
\begin{figure}
\begin{center}
\psfrag{pseudo_time_} {\scriptsize{Pseudo-time $t$}}
\psfrag{Dissipated_Energy_(J)} {\!\!\!\!\!\!\scriptsize{Dissipated energies (J/m)}}
\psfrag{Accumulated_Dissipation_due_to_damage} {\scriptsize{Accumulated dissipation due to damage}}
\psfrag{Accumulated_Dissipation_Pl} {\scriptsize{Accumulated dissipation}}
\psfrag{due_to_plasticity_case_(d)} {\scriptsize{due to plasticity, case (d)}}
\psfrag{(c)} {\!\!\!\!\!\!\!\!\scriptsize{cases (c),}}
\psfrag{(d)} {\scriptsize{ (d)}}
{\includegraphics[width=.7\textwidth, height=.50\textwidth]{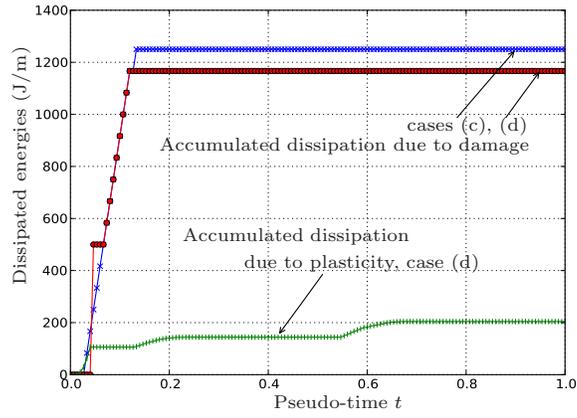}}
\caption{Evolution of the dissipated energies: (c)  interface damage and (d)  both interface damage and plasticity.}
\label{fig:ex1_5}
\end{center}
\end{figure}
\begin{figure}
\begin{center}
\psfrag{pseudo_time_} {\scriptsize{Pseudo-time $t$}}
\psfrag{Stored_Energy_(J)} {\!\!\!\!\!\!\!\!\!\!\!\!\scriptsize{Stored energies (J/m)}}
\psfrag{Interface_stored_energy_GII} {\tiny{Stored energy (tangent)}}
\psfrag{Interface_stored_energy_GI} {\tiny{Stored energy  (normal)}}
{\includegraphics[width=.7\textwidth, height=.50\textwidth]{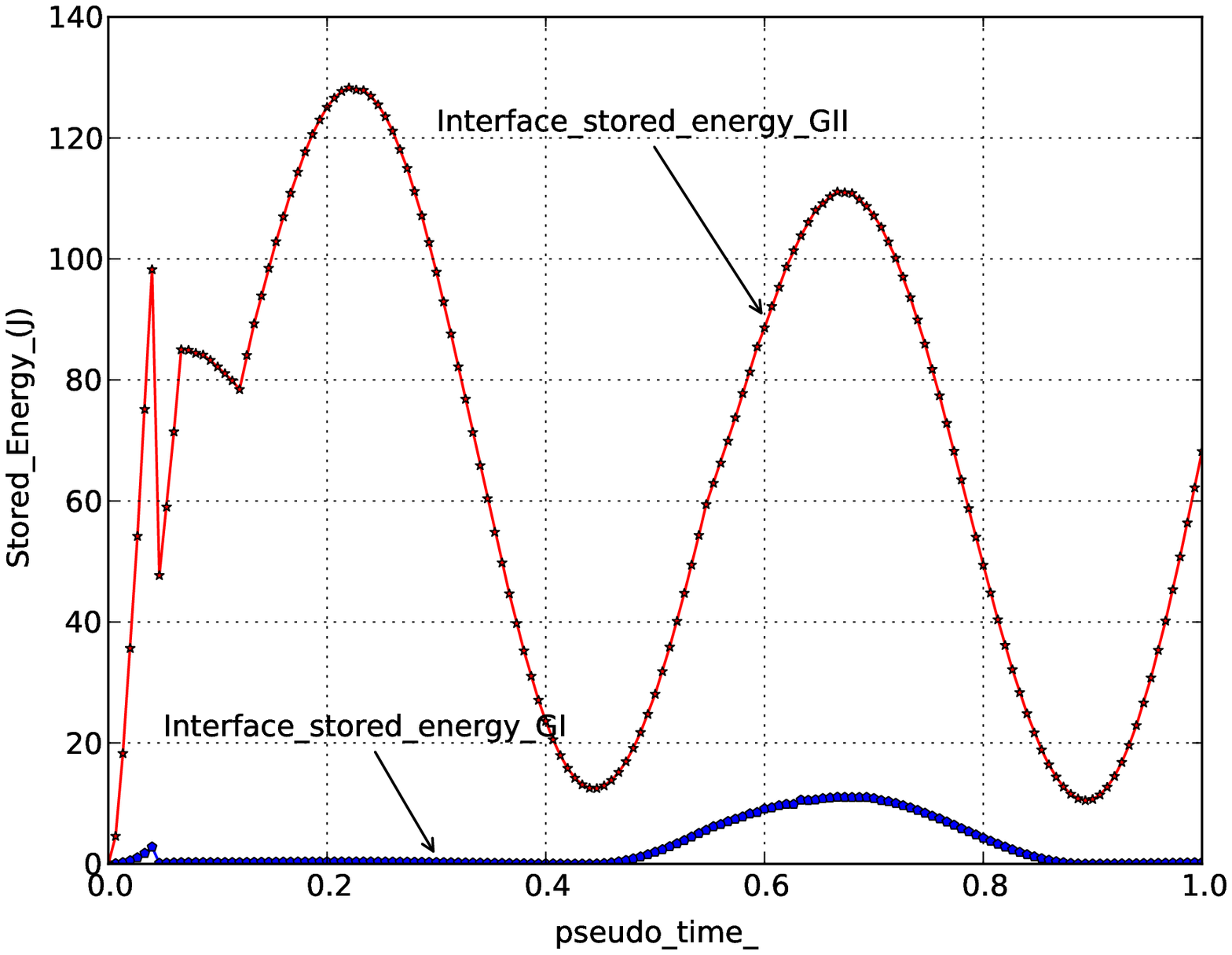}}
\caption{Evolution of the stored energies in the adhesive for case (d) considering both interface damage and plasticity.}
\label{fig:ex1_6}
\end{center}
\end{figure}
This behavior can be better understood from Fig.~\ref{fig:ex1_5}, where
the evolution of the accumulated dissipation with respect to the time $t$ is shown. In fact, $t$ is a kind of pseudo-time or process time
which  can arbitrarily be  re-scaled since the considered
system is rate-independent. Finally, the evolution of stored energies
in the adhesive layer due to opening and shear are shown in
Fig.~\ref{fig:ex1_6}.

\subsubsection{Traction instead of displacement loading}\label{par:TRACTIOMLOAD}
\begin{figure}
   \centering
   \def\svgwidth{.8\columnwidth}
   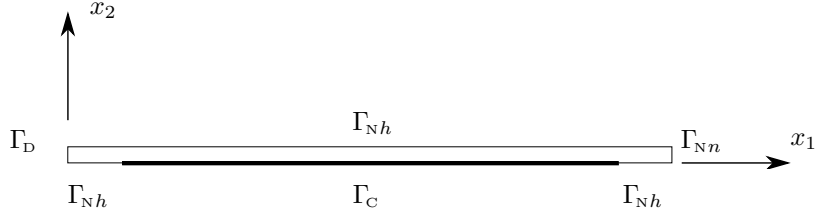
   \vspace*{2.5em}
   \caption{Modified geometry and boundary conditions used for the problem with prescribed non-zero tractions.}
\label{fig:modified}
\end{figure}
From an engineering point of view, we are highly interested in the case of external loading given by nonvanishing  Neumann   boundary conditions. This is because in numerous   experiments or   real applications, loading is described through external forces, moments  or tractions.
For these reasons    a modified  problem configuration shown in Fig.~\ref{fig:modified}  is studied in this section. The length of $\GC$ is $L_c=200$ mm, while the homogeneous Neumann boundary parts $\GN_h$ on the left  and right hand side of $\GC$ have lengths equal to $0.125 L_c$. The length and  height of $\Omega$ as well as the material properties are the same as in the previous example. The left vertical side of  $\Omega$ is fixed, defining $\GD$,  while uniform  normal tractions  applied on the right vertical side of  $\Omega$, defining $\GN_n$, are increasing in time, i.e. $p_1(t,x)=\psi(t) p_0$ and $p_2=0$ therein, with $p_0>0$ being a constant and $\psi(t)>0$   an increasing function, see Eq.\,\eqref{tfp}. Both interface damage and plasticity are considered.


\begin{figure}
\begin{center}
\psfrag{resultant_force_x___________________________(m)} {\!\!\!\!\!\!\scriptsize{Resultant horizontal force (KN/m)}}
\psfrag{displacement_(m)} {\!\!\!\!\!\!\!\!\!\!\!\!\!\!\!\!\!\!\!\!\!\!\!\!\scriptsize{Horizontal displacement (mm)}}
\psfrag{(L)} {\!\!\!\!\!\!\scriptsize{(L)}}
{\includegraphics[width=.7\textwidth, height=.50\textwidth]{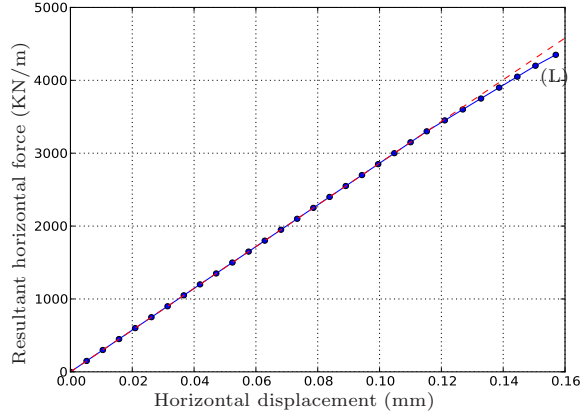}}
\caption{The horizontal resultant force on $\GN_n$ versus the horizontal  displacement of the  lower right corner of $\Omega$.}
\label{fig:ex2_1}
\end{center}
\end{figure}

The evolution of the horizontal resultant force on $\GN_n$ is plotted in Fig.~\ref{fig:ex2_1}.
 The computational analysis, including BTA together with the two-sided energy inequality checking,  stops at a point, marked in Fig.~\ref{fig:ex2_1} as point (L), where after some   plasticity development the first and at the same time the total damage of the adhesive layer occurs.   A dashed line in the same plot represents the tangent line  to the initial purely elastic part of the resultant force-displacement curve. Both lines separate due to the appearance of some plastic slip at $\GC$.  To capture a progressive damage propagation along $\GC$ would require   decreasing the applied load  after the peak  load is achieved.

\begin{figure}
\begin{center}
\psfrag{xxxxxxxxxxxxxxxxxxxxxxxxxxxx} {\tiny{According to Fig.~\ref{fig:springs}b}}
\psfrag{xxxxxxxxxxxxxxxxxxxxxxxxxxxy} {\tiny{10 elements on $\GC$}}
\psfrag{xxxxxxxxxxxxxxxxxxxxxxxxxxxw} {\tiny{50 elements on $\GC$}}
\psfrag{resultant_force_x___________________________(m)} {\,\,\,\scriptsize{Shear stress (MPa)}}
\psfrag{displacement_(m)} {\!\!\!\!\!\!\!\!\!\!\!\!\!\!\!\!\!\!\!\!\!\!\!\!\scriptsize{Tangential relative displacement (mm)}}
{\includegraphics[width=.7\textwidth, height=.50\textwidth]{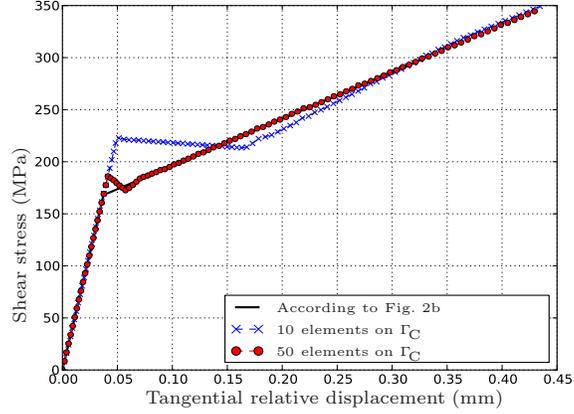}}
\caption{Shear stress-relative tangential displacement relation as computed at the rightmost point of the interface $\GC$. For this case, only plasticity is considered, damage being excluded.}
\label{fig:ex2_2}
\end{center}
\end{figure}

Finally, in order to analyse the pointwise behaviour of the numerical solution   at $\GC$, in particular regarding the evolution of plastic slip,  the same problem is solved again but including  interface plasticity only, i.e. no  damage at $\GC$  is possible. Numerically computed shear stress versus the  relative tangential displacement at the rightmost node of $\GC$ is shown in Fig.~\ref{fig:ex2_2}, together with the  expected tangential stress-relative displacement law of Fig.~\ref{fig:springs}. An ``\emph{overshooting}'' phenomenon takes place when plasticity occurs, a similar behaviour may also be observed in Fig.~\ref{fig:ex1_2}. This phenomenon is essentially associated to the time and spatial  discretization of the problem, in particular  possibly due to some   oscillations of the traction solution near the crack tip, and therefore  can   gradually  be eliminated by a spatial discretization refinement as observed in Fig.~\ref{fig:ex2_2}.

\subsection{Practical application}\label{par:APPLICATION}
 The problem configuration shown in Fig.~\ref{fig_m1} is considered again. A monotonic hard-device loading is assumed by prescribing horizontal and vertical displacements, respectively, as $u_1(t,x)=t\,w_1(x)$ with $w_1(x)=0.6$mm and $u_2(t,x)=0.6u_1(t,x)$  at the  right-hand side  of the rectangle $\Gamma$, defining the Dirichlet boundary $\GD$. All the other boundary parts are considered to be traction free, defining the Neumann boundary $\GN$, except for the contact surface $\GC$. The problem evolution is  represented as function of  the fictitious time $t$.

 An advantage of the present method considering a continuous distribution of springs along the interface in comparison with the classical  fracture mechanics, which assumes a perfect interface except for some cracked zones, is that no special mesh refinement is needed near the crack tip and a uniform mesh can be employed along $\GC$, similarly as in the Cohesive Zone Models, see \cite{HuRuLoJa2011JA} and further references therein. The tractions along the adhesive layer of the present type are bounded although  traction concentrations can be expected at the end-points of the adhesive layer, which may correspond to crack tips, cf. \cite{LENCI2001,TMGCP10ACTA}.
 Actually, in the present model, these tractions are limited by  the critical values of normal and tangential tractions $\sigma_{n,{\rm crit}}$ and $\sigma_{t,{\rm crit}}$, respectively. If the  adhesive layer at $\GC$ in the present problem in Fig.~\ref{fig_m1}  would be replaced by perfect bonding conditions, stress singularities would appear at both extremes of the bonded part $\GC$. The stress singularity at the right extreme of  $\GC$,   corresponding to the classical oscillatory singularity of the open model of  an  interface crack   between an elastic and an infinitely rigid solid, cf. \cite{Williams1959}, would be more severe than that at the left extreme. In such a case, a strongly refined mesh or special singular elements would be needed for a proper problem discretization of the crack tip neighbourhood. It should be mentioned here that intuitive refinement without having rigorous local error indicators is sometimes dangerous and may destroy convergence which is standardly guaranteed on uniformly refined meshes only.
 
\begin{figure}
\begin{center}
\psfrag{xxxxxxxxxx} {\scriptsize{$n$=54}}
\psfrag{xxxxxxxxxy} {\scriptsize{$n$=108}}
\psfrag{xxxxxxxxxz} {\scriptsize{$n$=216}}
\psfrag{resultant_force_x___________________________(m)} {\scriptsize{Normal tractions $t_2$ (MPa)}}
\psfrag{resultant_force_y___________________________(m)} {\scriptsize{Tangential tractions $t_1$ (MPa)}}
\psfrag{displacement_(m)} {\!\!\!\!\!\!\!\!\!\!\!\!\!\!\!\!\scriptsize{Coordinate $x$ on $\GD$ (mm)}}
{\includegraphics[width=.7\textwidth, height=.50\textwidth]{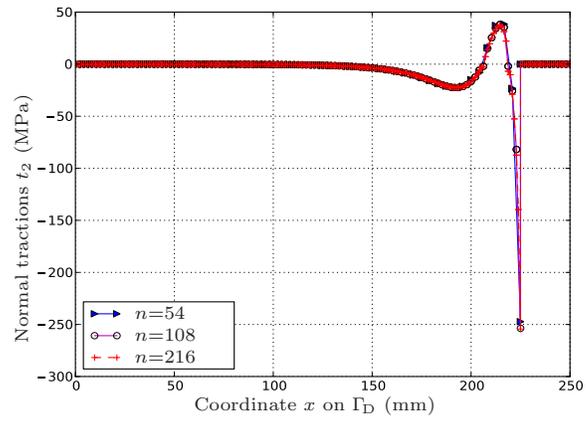}} (a)
{\includegraphics[width=.7\textwidth, height=.50\textwidth]{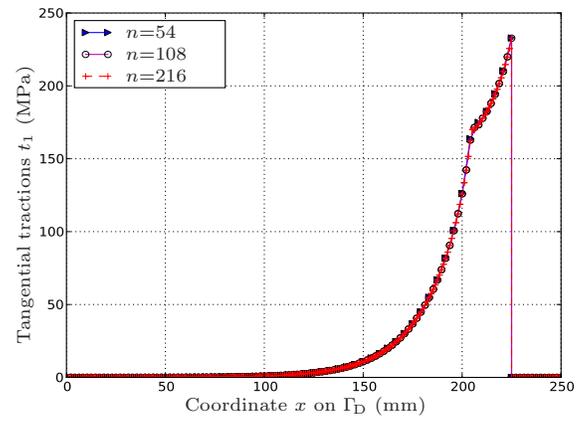}} (b)
\caption{(a) Normal and (b) tangent tractions along the adhesive zone $\GC$, for the three finest uniform boundary element meshes, $n$ the number of boundary elements at $\GC$.}
\label{fig:ex3_m_1}
\end{center}
\end{figure}

\begin{figure}
\begin{center}
\psfrag{xxxxxxxxxxxxxxxxxxxxxxxxxxxxxxxxxxxxxxxxxw} {\scriptsize{(a) Strain energy of the bulk}}
\psfrag{xxxxxxxxxxxxxxxxxxxxxxxxxxxxxxxxxxxxxxxxxx} {\scriptsize{(b) Total energy}}
\psfrag{xxxxxxxxxxxxxxxxxxxxxxxxxxxxxxxxxxxxxxxxxy} {\scriptsize{(c) Maximum absolute normal traction}}
\psfrag{percentage_error} {\scriptsize{Percentage error}}
\psfrag{number_of_elements} {\scriptsize{Number of elements on $\GC$}}
\psfrag{displacement_(m)} {\!\!\!\!\!\!\!\!\!\!\!\!\!\!\!\!\scriptsize{Coordinate $x$ on $\GD$ (mm)}}
{\includegraphics[width=.7\textwidth, height=.50\textwidth]{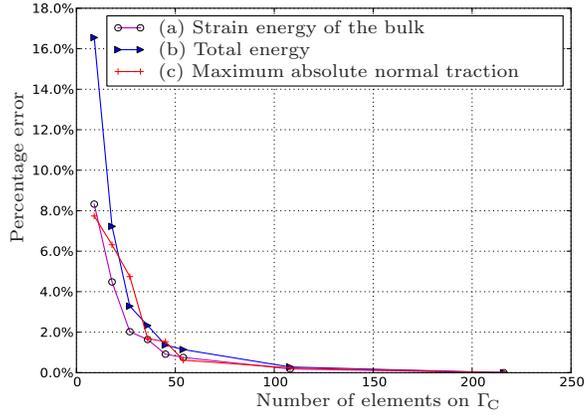}}
\caption{Percentage difference of (a) the total energy, (b) the elastic strain energy and (c) the maximum normal tractions at the adhesive zone $\GC$, for  uniform, along the horizontal and vertical edges, boundary element meshes.}
\label{fig:ex3_m_2}
\end{center}
\end{figure}

 Nevertheless, in the present model, the fact that a local mesh refinement  at the crack tip is not needed makes easy the modeling of damage progression with the  crack tip moving along the interface. 
 In order to check this statement, we have solved the problem in Fig.~\ref{fig_m1}  by using a series of uniform boundary element meshes, the three finest meshes having 126 (60 and 3 elements along each horizontal  and vertical side, respectively), 252 and 504   elements on $\Gamma$ (that corresponds to 54, 108 and 216 elements on $\GC$).  The traction solutions along $\GC$ for these three finest meshes shown in Fig.~\ref{fig:ex3_m_1}  correspond to horizontal prescribed displacement $u_1=0.28$mm, when no damage appears although some interface plasticity has evolved. A strong traction concentration at the right extreme of $\GC$ can be observed in these plots which, however, indicate that even for the coarsest mesh (54 elements on $\GC$)  the solution obtained is sufficiently accurate for the purpose of the present study. An additional checking is shown  in Fig.~\ref{fig:ex3_m_2}, where the percentage differences of the computed strain energy in the bulk,  the computed total energy (that is the sum of the stored energy  and the dissipated energy at time $t$, $\calE(t,u(t),z(t))+\Diss_{\mathcal R}(z;[0,t]$, see \eqref{total-energy}), and the maximum absolute value of  normal tractions at $\GC$ computed  by a coarse mesh and the finest mesh (216 elements on $\GC$) are plotted. These plots confirm  that the percentage difference of the strain energy and of the maximum normal traction for the mesh with 54 elements on $\GC$ is sufficiently small, in particular it is about 1\%. Therefore, this mesh is used in the following complete numerical study of the present problem.
%
\begin{figure}
\begin{center}
\psfrag{pseudo_time_} {\scriptsize{Pseudo-time $t$}}
\psfrag{Energy_(J)} {\!\!\!\!\!\!\!\!\!\!\!\!\scriptsize{Energies (J/m)}}
\psfrag{Bulk} {\tiny{Bulk}}
\psfrag{Interface} {\tiny{Interface}}
\psfrag{Dissipated} {\tiny{Dissipated}}
\psfrag{Lower estimate} {\tiny{Lower estimate}}
\psfrag{Total} {\tiny{Total}}
\psfrag{Upper estimate} {\tiny{Upper estimate}}
{\includegraphics[width=.7\textwidth, height=.50\textwidth]{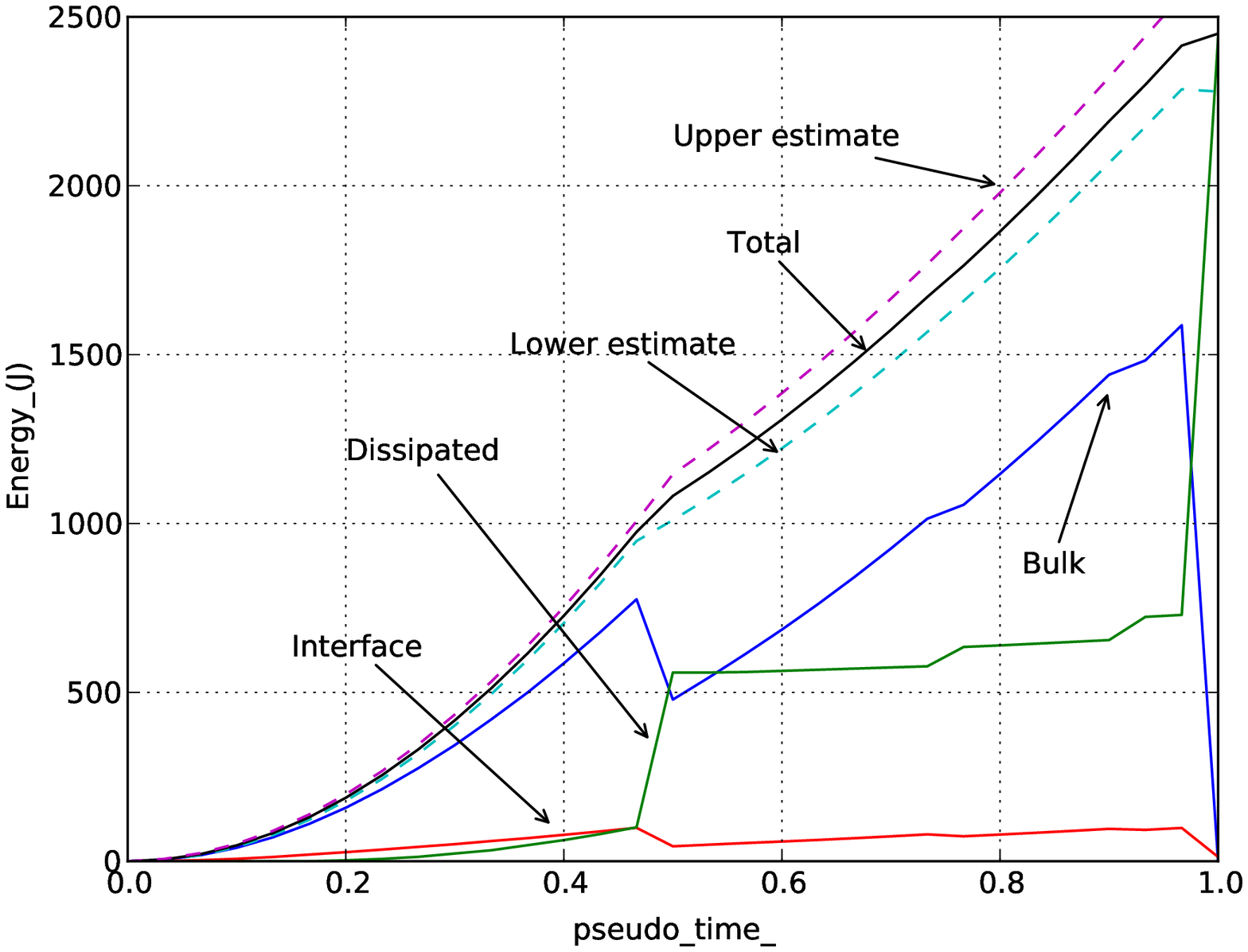}}
\caption{Evolution of the energies.}
\label{fig:ex3_1}
\end{center}
\end{figure}

\begin{figure}
\begin{center}
\psfrag{displacementX_(m)} {\scriptsize{Coordinate $x_1$ (m)}}
\psfrag{displacementY_(m)} {\scriptsize{Coordinate $x_2$ (m)}}
\psfrag{xxxxxxxxxxxxxxxxxxxxxxxxxxxx} {\tiny{Undeformed shape}}
\psfrag{xxxxxxxxxxxxxxxxxxxxxxxxxxxy} {\tiny{Deformed shape before damage}}
\psfrag{xxxxxxxxxxxxxxxxxxxxxxxxxxxw} {\tiny{Deformed shape after damage}}
{\includegraphics[width=.7\textwidth, height=.50\textwidth]{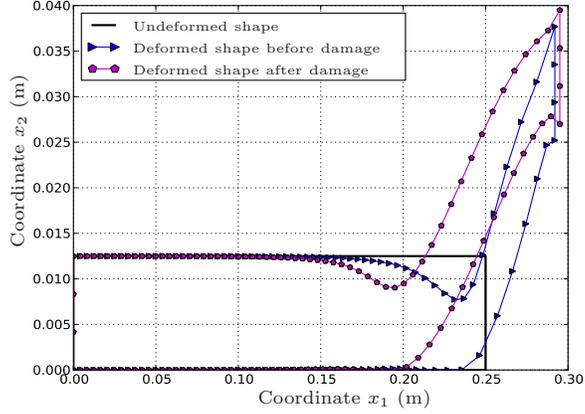}}
\caption{Deformed shape of the elastic domain for two time steps, just before and after damage initiation. A scale factor of 150 has been used.}
\label{fig:ex3_2}
\end{center}
\end{figure}

Fig. \ref{fig:ex3_1} shows the evolution of different energies  computed, in particular, the energy stored in the elastic bulk,  in the adhesive layer and the dissipated energy. Also  the total energy, which is actually minimized in the time stepping procedure, together with the lower and upper estimates of energy are shown. As it can be seen in Fig. \ref{fig:ex3_1}, the global minimization procedure defines   the end   of the delamination process, where the remaining undamaged part of the adhesive layer is debonded, as the point where the sum of the stored energy  in the bulk and the adhesive layer   is basically equal to the energy needed to   delaminate the   undamaged part of the adhesive layer, which is given by the dissipated energy in the last time step. In Fig. \ref{fig:ex3_2}, the deformed shape of the bulk, is plotted for two time steps, just before and after the first damage, respectively.

\begin{figure}
\begin{center}
\psfrag{xxxxxxxxxx} {\tiny{BEM}}
\psfrag{xxxxxxxxxy} {\tiny{Interface}}
\psfrag{resultant_force_x___________________________(m)} {\scriptsize{Normal tractions $t_2$ (MPa)}}
\psfrag{displacement_(m)} {\!\!\!\!\!\!\!\!\!\!\!\!\!\!\!\!\scriptsize{Coordinate $x_1$ on $\GC$ (mm)}}
{\includegraphics[width=.7\textwidth, height=.50\textwidth]{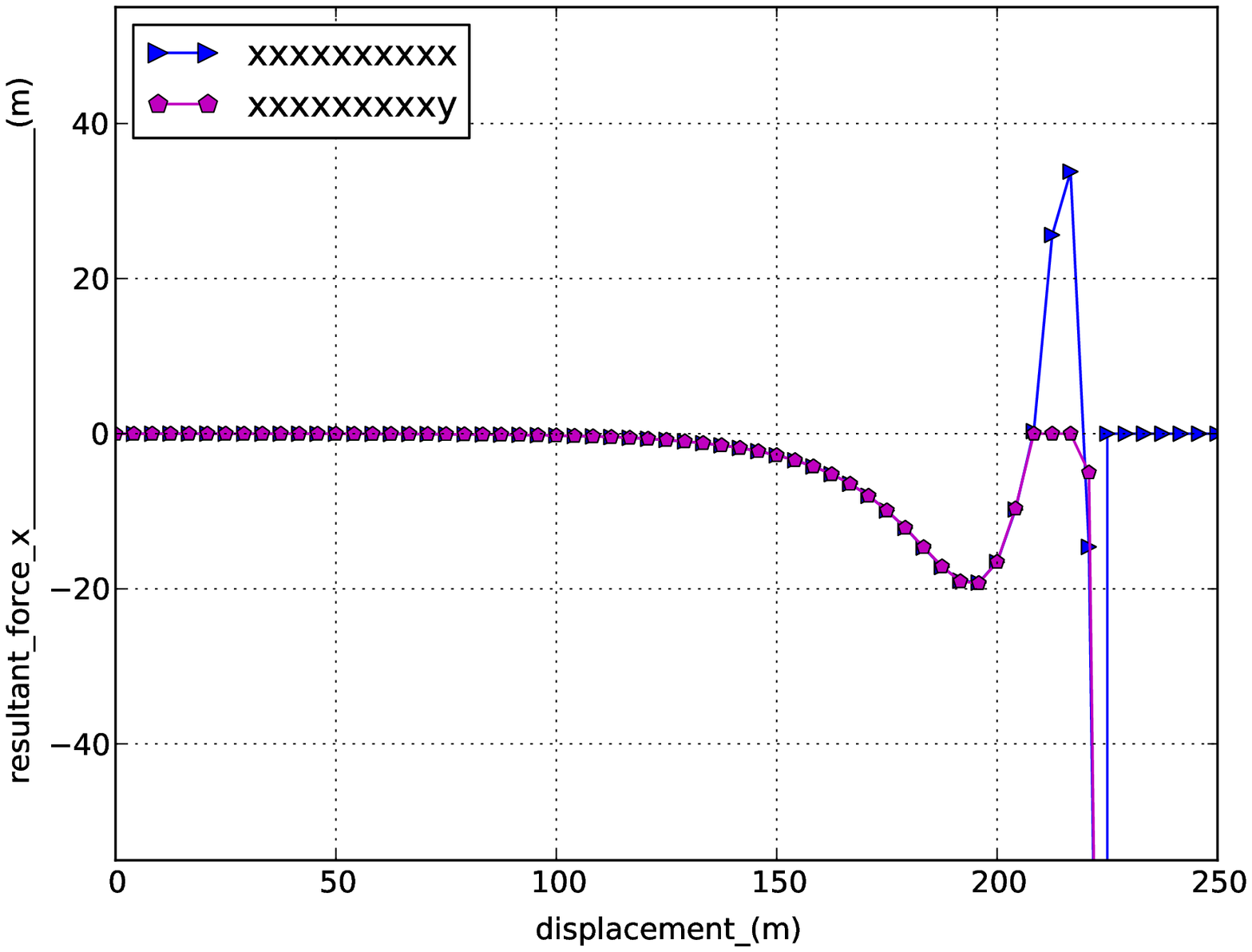}} (a)
{\includegraphics[width=.7\textwidth, height=.50\textwidth]{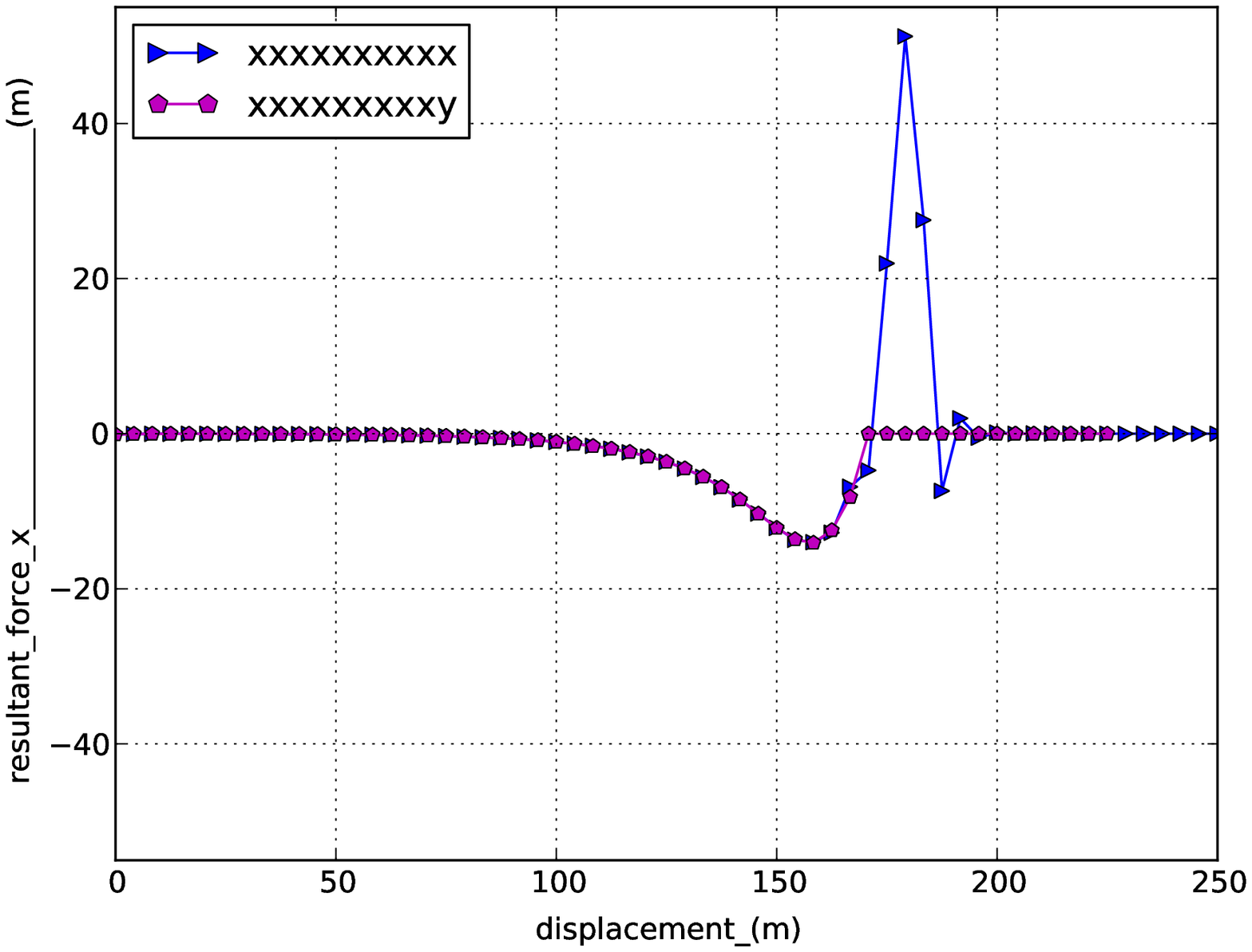}} (b)
\caption{Normal tractions along the adhesive zone $\GC$, just (a) before and (b) after the first crack opening, computed by BEM as well as in the interface elements.}
\label{fig:ex3_3}
\end{center}
\end{figure}

\begin{figure}
\begin{center}
\psfrag{xxxxxxxxxx} {\tiny{BEM}}
\psfrag{xxxxxxxxxy} {\tiny{Interface}}
\psfrag{resultant_force_x___________________________(m)} {\scriptsize{Tangential tractions $t_1$ (MPa)}}
\psfrag{displacement_(m)} {\!\!\!\!\!\!\!\!\!\!\!\!\!\!\!\!\scriptsize{Coordinate $x_1$ on $\GC$ (mm)}}
{\includegraphics[width=.7\textwidth, height=.50\textwidth]{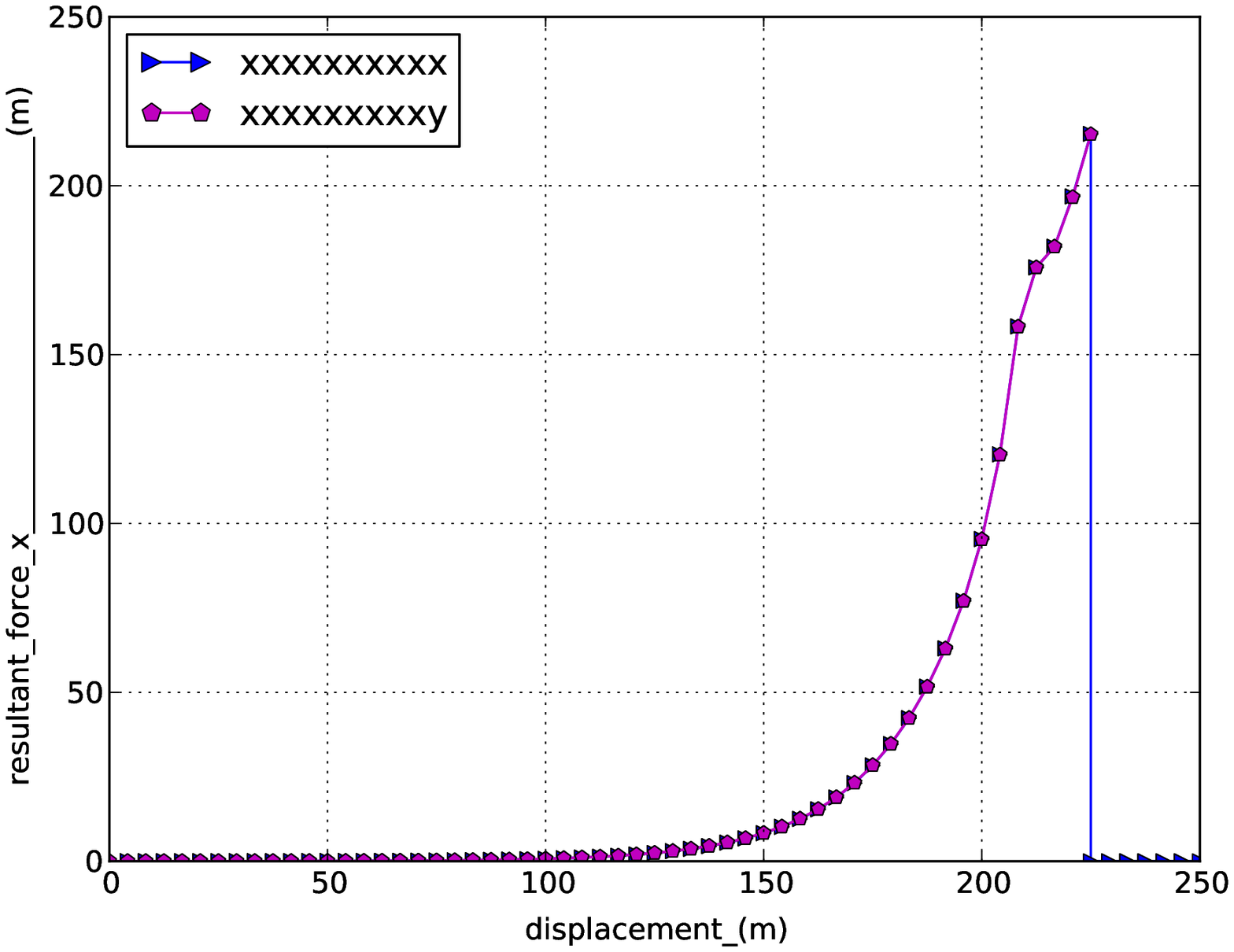}} (a)
{\includegraphics[width=.7\textwidth, height=.50\textwidth]{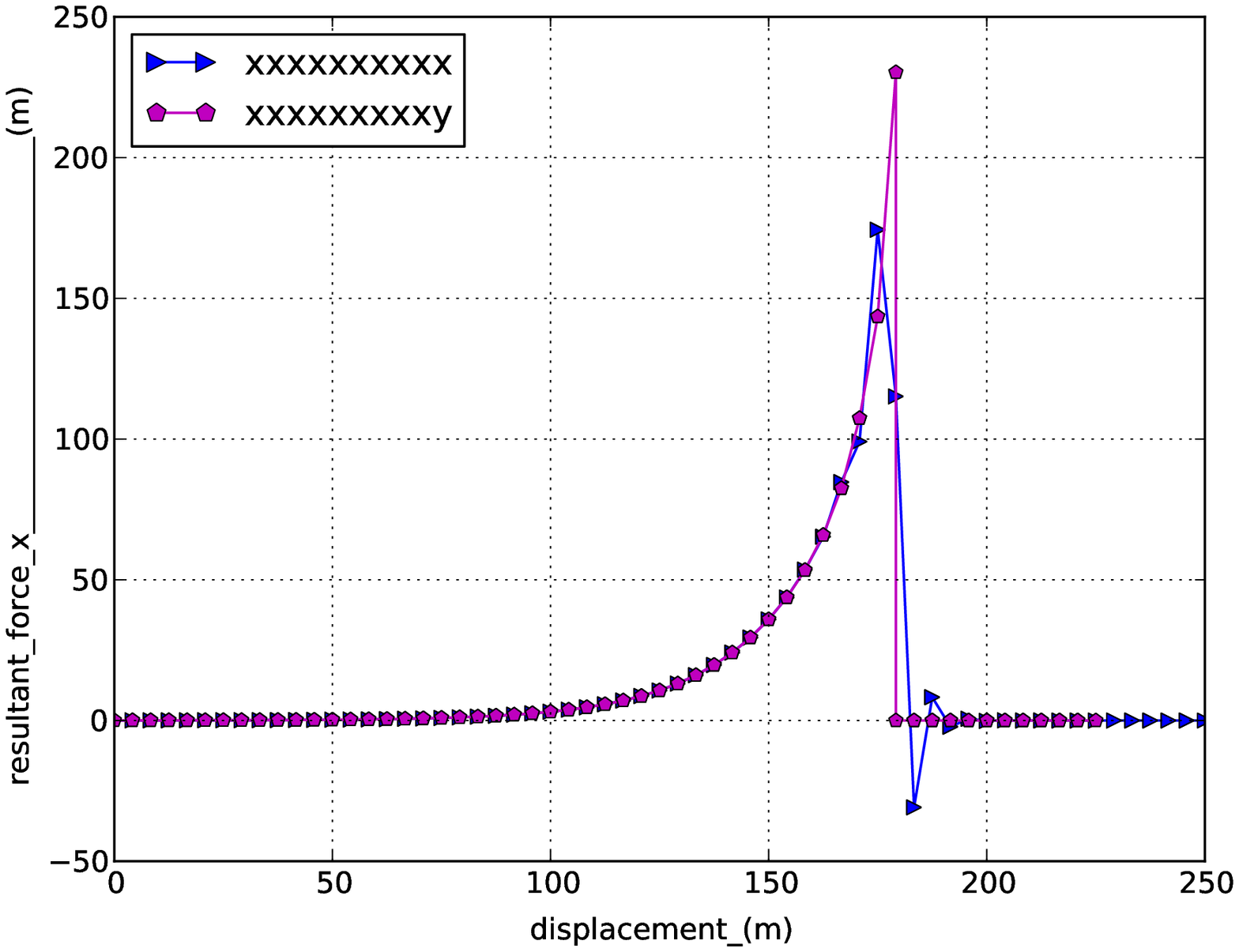}} (b)
\caption{Tangential tractions along the adhesive zone $\GC$, just (a) before and (b) after the first crack opening, computed by BEM as well as in the interface elements.}
\label{fig:ex3_4}
\end{center}
\end{figure}

Figs. \ref{fig:ex3_3} and \ref{fig:ex3_4} present the components of normal and tangential traction  vector along  the adhesive zone $\GC$. A very good agreement exists between the computed tractions for each subdomain by BEM   and those computed in the adhesive layer, although the equilibrium has not been imposed directly but it results as a consequence of the energy minimization. Progressive extension of the traction free portion of the original $\GC$ because of the damage propagation ($\zeta=0$) can be observed in Figs. \ref{fig:ex3_3}b) and \ref{fig:ex3_4}b). It should be mentioned   that the  portion  of $\GC$  which is totally damaged is still kept as a part of the minimization procedure, where nodal displacement values participate as  unknowns in the minimization procedure and their values are used in the BEM solution of the pertinent BVP. For this reason, in fact      an approximation of the developed traction-free    zone is computed by BEM for each subdomain. Obviously other algorithms might be used where after the total damage of a portion of the adhesive layer a switch in the type of the boundary condition (e.g. from Dirichlet to vanishing Neumann boundary condition) is taken into account along this boundary portion in the BEM computation for each subdomain. Nevertheless, we have been interested in the results obtained by the present simple   procedure.  Normal compressive tractions computed by BEM can be observed in  Fig. \ref{fig:ex3_3} in zones where  vanishing normal tractions are obtained in the interface elements. This is due to the fact that the rigid obstacle undertakes these compressive tractions.

\begin{figure}
\begin{center}
\psfrag{xxxxxxxxxxxxxxxxxxxx} {\tiny{Signorini contact}}
\psfrag{xxxxxxxxxxxxxxxxxxxy} {\tiny{Damage and plasticity}}
\psfrag{resultant_force_y___________________________(m)} {\!\!\!\!\!\!\scriptsize{Resultant vertical force (KN/m)}}
\psfrag{resultant_force_x___________________________(m)} {\!\!\!\!\!\!\scriptsize{Resultant horizontal force (KN/m)}}
\psfrag{displacement_(m)} {\scriptsize{Pseudo-time $t$}}
\psfrag{(A)} {\tiny{(A)}}
\psfrag{(B)} {\tiny{(B)}}
\psfrag{(C)} {\tiny{(C)}}
\psfrag{(D)} {\tiny{(D)}}
\psfrag{(E)} {\tiny{(E)}}
{\includegraphics[width=.7\textwidth, height=.50\textwidth]{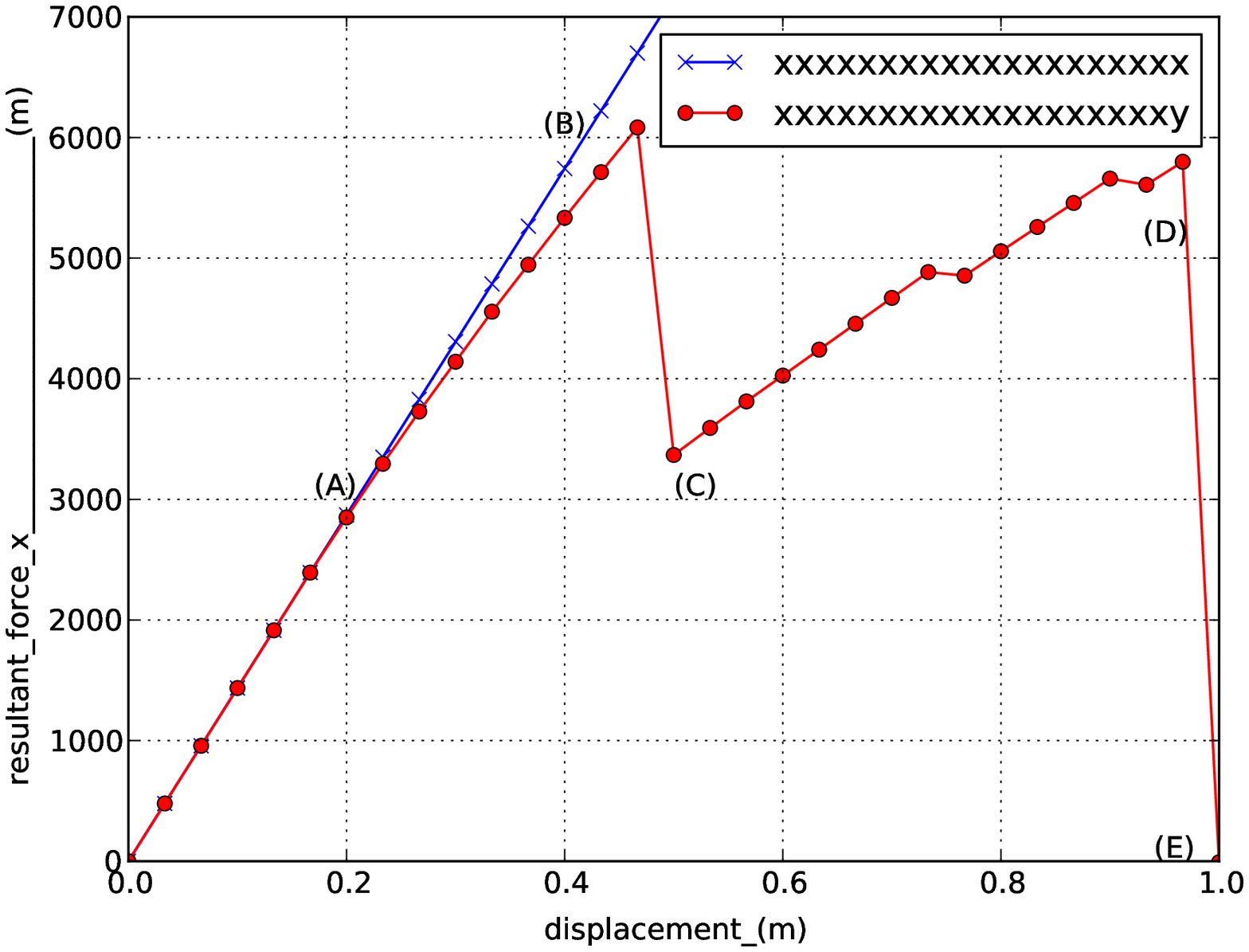}} (a)
{\includegraphics[width=.7\textwidth, height=.50\textwidth]{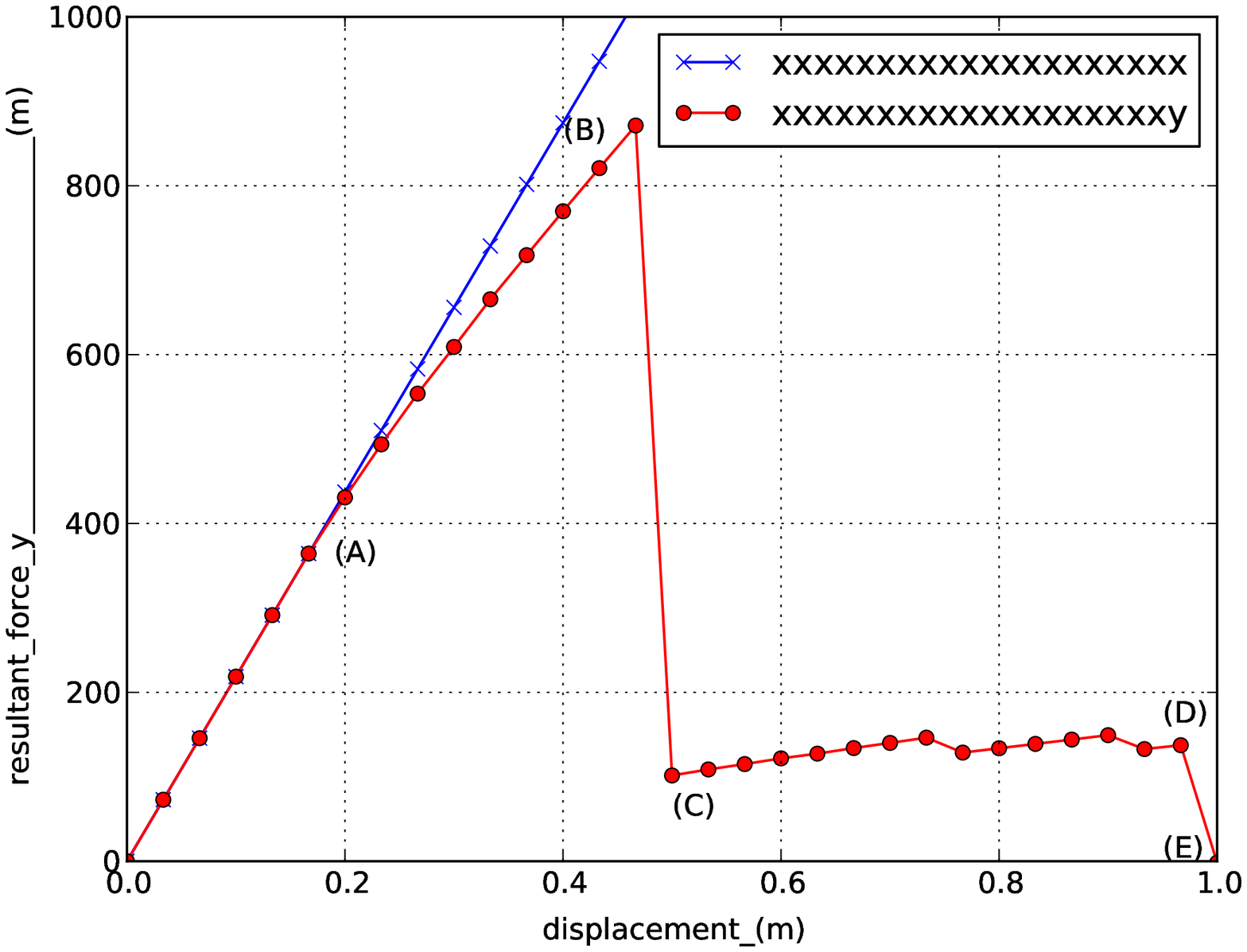}} (b)
\caption{Evolution of (a) the horizontal and (b) vertical resultant force on $\GD$.}
\label{fig:ex3_5}
\end{center}
\end{figure}

Finally, in Fig. \ref{fig:ex3_5} the resultant forces acting at $\GD$ versus the  pseudo-time $t$ are shown. These two plots have some similarities in the behaviour that may come up through the characteristic points (A)-(E). Up to point (A) the linear elastic behaviour manifests in both the solid and  adhesive, at this point a plastic slip appears. Then, at point (B) the first damage appears in the first 11 boundary elements that are situated on the right-hand side of the adhesive layer. This new crack length results in a ``jump down'' of the resultant forces up to point (C). Then, up to point (D) the damaged zone is progressively extended and finally after point (D) the remaining adhesive zone is   damaged instantaneously. The problem evolution ends up at point (E), where a rigid body motion of the elastic body takes place. The increment of the crack length  from point (B) to (C)  equals   45.83 mm. In the same figures also the linear   elastic responses, obtained for the same configuration taking into account only Signorini contact without any interface damage and plasticity, are plotted in order to facilitate the observation of  the  initiation of plasticity and/or damage.

\section{Conclusions}
A boundary element implementation of a computational procedure based on
\COL{an} energetic\COL{-}solution framework for the delamination problems has been presented.
A specific model for the adhesive interfaces, which distinguishes the amount of energy
dissipated in opening Mode I and shear Mode II has been adopted. This model
involves two inelastic  internal  variables  on delaminating surfaces,
namely the damage variable $\zeta$ and the plastic slip $\pi$.
Some details regarding the formulation of the collocation BEM as well as
the optimization procedures necessary for solving the global
minimization problem, inherent in the formulation, have been discussed.
A few numerical tests  have been presented in order to
analyse the behaviour of the present delamination model and performance of the algorithms implemented in a collocational BEM code.





\newpage
\bibliographystyle{spmpsci}      
\bibliography{PMR_IABEM_CM}







\end{document}